\documentclass{amsart}

\usepackage{amssymb}
\usepackage{amscd}
\usepackage{amsmath}
\usepackage{latexsym}
\usepackage{verbatim}
\usepackage[latin1]{inputenc}
\usepackage[dvips]{graphics}

\relax
\usepackage[all]{xy}

\newdir{ >}{{}*!/-10pt/\dir{>}}

{\begin{list}{}{%
\settowidth{\labelwidth}{\textsf{#1}}%
\setlength{\leftmargin}{\labelwidth}
\addtolength{\leftmargin}{0pt}}}%
{\end{list}}

\newtheorem{Thm}{Theorem}[section]
\newtheorem{Prop}[Thm]{Proposition}
\newtheorem{PropDef}[Thm]{Proposition-Definition}
\newtheorem{Lem}[Thm]{Lemma}
\newtheorem{Cor}[Thm]{Corollary}
\newtheorem{CorDef}[Thm]{Corollary-Definition}

\theoremstyle{definition}
\newtheorem{Def}[Thm]{Definition}
\newtheorem{Not}[Thm]{Notation}

\newtheorem{Ex}[Thm]{Example}

\newtheorem{Conv}[Thm]{Convention}

\theoremstyle{remark}
\newtheorem{Rem}[Thm]{Remark}

\newcommand{\Prf}{\noindent\textit{Proof. }}
\newcommand{\PrfOf}[1]{\noindent\textit{Proof of #1.}}

\newcommand{\bbN}{\mathbb{N}}

\newcommand{\bbE}{\mathbb{E}}
\newcommand{\bbF}{\mathbb{F}}

\newcommand{\mor}{\mathop{\rm mor}\nolimits}
\newcommand{\obj}{\mathop{\rm obj}\nolimits}

\newcommand{\res}{\mathop{\rm res}\nolimits}
\newcommand{\ind}{\mathop{\rm ind}\nolimits}

\newcommand{\incl}{\mathop{\rm incl}\nolimits}
\newcommand{\colim}{\mathop{\rm colim}}
\newcommand{\Lcolim}{\mathop{L\rm colim}}
\newcommand{\hocolim}{\mathop{\rm hocolim}}
\newcommand{\id}{\mathop{\rm i\hspace*{-.03em}d}\nolimits}

\newcommand{\arr}{\mathop{\rm arr}\nolimits}
\newcommand{\too}{\longrightarrow}

\newcommand{\nspace}{\hspace*{-.1em}}
\newcommand{\nnspace}{\hspace*{-.05em}}
\newcommand{\nnnspace}{\hspace*{-.02em}}

\newcommand{\pspace}{\hspace*{.05em}}

\newcommand{\cat}[1]{{\mathcal{#1}}}

\newcommand{\smallbullet}{{\scriptscriptstyle \bullet}}

\newcommand{\Top}{\mathcal{T}\!\nspace{\rm o\nnspace p}}

\newcommand{\Sets}{\mathcal{S}\nspace{\rm e\nnspace t\nnspace s}}
\newcommand{\sSets}{{\rm s}\Sets}

\newcommand{\DDelta}{\mathbf{\Delta}}
\newcommand{\op}{^{\nnnspace\rm o\nnspace p}}

\newcommand\UUU[3]{\cat{U}_{#1}(#2,#3)}

\newcommand\USC[2]{\cat{U}_{\cat #1}(\cat #2)}

\newcommand\UUSC[2]{\cat{U}_{#1}(#2)}

\newcommand\USCD[3]{\UUU{\cat #1}{\cat #2}{\cat #3}}

\newcommand\FDC[2]{\mor_{\cat #1,\cat #2}} 
\newcommand\FFDC[2]{{\bbF}_{\cat #1}^{\,\cat #2}}
\newcommand\FFC[1]{{\bbF}\kern.1em{\cat #1}}

\newcommand{\ie}{{\sl i.e.\ }}
\newcommand{\eg}{{\sl e.g.\ }}

\newcommand{\weqs}{weak\ equivalences}

\newcommand{\itemspace}{\qquad\quad\;\;}

\newcommand{\coeq}{\mathop{\rm coeq}\nolimits}

\newcommand{\Map}{\mathop{\rm Map}\nolimits}

\newcommand{\Cst}[2]{{\rm C}\nnspace{\rm s}\nnspace{\rm t}_{#1}^{#2}}

\newcommand{\cterminal}{c_{\nspace{\scriptscriptstyle\infty}}}

\newcommand{\stacksym}[2]{\renewcommand{\arraystretch}{0.1}
\begin{array}{@{}c@{}}
#2\\ #1\\
\end{array}
\renewcommand{\arraystretch}{1}}

\newcommand{\adjtoo}{\ \stacksym{\longleftarrow}{\longrightarrow}\ }

\newcommand{\noloc}{\hspace*{.05em}:\hspace*{-.13em}}

\newcommand{\rmap}[1]{\mathop{\rm map}_{\nspace{\scriptscriptstyle#1}}\nolimits}
\newcommand{\lmap}[1]{{}_{{\scriptscriptstyle#1}}\!\nnspace\mathop{{\rm map}}\nolimits}

\newcommand{\stackup}[2]{\renewcommand{\arraystretch}{0.4}
\begin{array}[b]{@{}c@{}}
{\scriptscriptstyle #1}\\ #2\\
\end{array}\,
\renewcommand{\arraystretch}{1}}

\newcommand{\stackdown}[2]{\renewcommand{\arraystretch}{0.4}
\begin{array}[t]{@{}c@{}}
  #2\\[.06em] {\scriptscriptstyle #1}\\
\end{array}\,
\renewcommand{\arraystretch}{1}}

\newcommand{\stackthree}[3]{\renewcommand{\arraystretch}{0.4}
\begin{array}{@{}c@{}}
  {\scriptscriptstyle #1}\\ #3\\[.06em] {\scriptscriptstyle #2}\\
\end{array}\,
\renewcommand{\arraystretch}{1}}

\newcommand{\uCdot}[2]{\stackup{\cat{#1}}{#2}}
\newcommand{\dCdot}[2]{\stackdown{\cat{#1}}{#2}}
\newcommand{\udCDdot}[3]{\,\stackthree{\cat{#1}}{\cat{#2}}{#3}}

\newcommand{\smallstackup}[2]{\renewcommand{\arraystretch}{0.2}
\begin{array}[b]{@{}c@{}}
{\scriptscriptstyle #1}\\ {\scriptscriptstyle #2}\\
\end{array}
\renewcommand{\arraystretch}{1}}

\newcommand{\smallstackdown}[2]{\renewcommand{\arraystretch}{0.2}
\begin{array}[t]{@{}c@{}}
  {\scriptscriptstyle #2}\\[.08em] {\scriptscriptstyle #1}\\
\end{array}
\renewcommand{\arraystretch}{1}}

\newcommand{\smallstackthree}[3]{\renewcommand{\arraystretch}{0.2}
\begin{array}{@{}c@{}}
  {\scriptscriptstyle #1}\\ {\scriptscriptstyle #3}\\[.08em] {\scriptscriptstyle #2}\\
\end{array}
\renewcommand{\arraystretch}{1}}

\newcommand{\smalluCdot}[2]{\pspace\smallstackup{\cat{#1}}{#2}}
\newcommand{\smalldCdot}[2]{\pspace\smallstackdown{\cat{#1}}{#2}}
\newcommand{\smalludCDdot}[3]{\pspace\smallstackthree{\cat{#1}}{\cat{#2}}{\cat{#3}}}

\newcommand{\pair}[2]{\hspace*{.075em}#1\;#2}

\newcommand\mQC[1]{{\mathcal Q}_{\cat #1}}
\newcommand\mQCD[2]{{\mathcal Q}_{\cat #2}^{\cat #1}}
\newcommand\nQX[1]{{\mathcal Q}#1}
\newcommand\xiC[1]{\xi^{\cat #1}}

\newcommand\xiCX[2]{\xi_{#2}^{\cat #1}}
\newcommand\xiCD[2]{\xi^{\cat #1,\cat #2}}
\newcommand\xiCDX[3]{\xiCD{#1}{#2}_{#3}}
\newcommand\nzetaX[1]{\zeta_{#1}}

\newcommand\deco[1]{{}^{{}^{\natural}}\kern-0.1em #1} 

\newcommand\Qbar{\smash[t]{\overline{Q}}}
\newcommand\QCbar[1]{\Qbar_{\nspace\cat #1}}
\newcommand\QCbarpr[1]{\deco{\Qbar}_{\nspace\cat #1}}
\newcommand\QDCbar[2]{\Qbar_{\cat #2}^{\,\cat #1}}
\newcommand\QDCbarpr[2]{\deco{\Qbar}_{\cat #2}^{\,\cat #1}}
\newcommand\xiCbar[1]{\bar{\xi}^{\cat #1}}
\newcommand\xiCbarpr[1]{\deco{\bar{\xi}}^{\cat #1}}
\newcommand\xiCXbar[2]{\bar{\xi}_{#2}^{\cat #1}}

\newcommand\xiCDXbar[3]{\bar{\xi}_{#3}^{\cat #1,\cat #2}}
\newcommand\xiCDXbarpr[3]{\deco{\bar{\xi}}_{#3}^{\cat #1,\cat #2}}

\newcommand\EEC[1]{{\bbE}\cat #1}
\newcommand\EECpr[1]{\deco{\bbE}\cat #1}
\newcommand\EEDC[2]{{\bbE}_{\cat #1}^{\cat #2}}
\newcommand\EEDCpr[2]{\deco{\bbE}_{\cat #1}^{\cat #2}}
\newcommand\vthetapr{\deco{\kern-0.1em\vartheta}}
\newcommand\vthetaC[1]{\vartheta_{\cat #1}}
\newcommand\vthetaDC[2]{\vartheta_{\cat #1}^{\cat #2}}
\newcommand\vthetaCpr[1]{\vthetapr_{\cat #1}}
\newcommand\vthetaDCpr[2]{\vthetapr_{\cat #1}^{\cat #2}}

\newcommand{\comma}{\!\searrow\!}
\newcommand{\smallcomma}{\searrow}
\newcommand{\commaC}[1]{\!\!{\renewcommand{\arraystretch}{0.2}
\begin{array}[b]{c}
{\scriptscriptstyle \;\;\cat #1} \\[-.4em]
\searrow \\
\end{array}
\renewcommand{\arraystretch}{1}}\!\!}
\newcommand{\smallcommaC}[1]{\!\!\!\nspace{\renewcommand{\arraystretch}{0.2}
\begin{array}[b]{c}
{\scriptscriptstyle \;\;\cat #1} \\[-.3em]
{\scriptstyle \searrow} \\
\end{array}
\renewcommand{\arraystretch}{1}}\!\!\!\nspace}

\newcommand{\qquest}{?\hspace*{-.06em}?}
\newcommand{\qqquest}{?\hspace*{-.06em}?\hspace*{-.06em}?}
\newcommand{\smallqquest}{?\hspace*{-.08em}?}

\newcommand{\sk}{\mathop{\rm sk}\nolimits}
\newcommand{\nd}{\mathop{\rm nd}\nolimits}

\newcommand{\Projec}{\rm P\nnspace r}

\newcommand{\coeqcat}{\xymatrix@C=1.2em{\smallbullet\!\ar@<-.22em>[r] \ar@<.22em>[r]&\!\smallbullet}}

\newcommand{\simtoo}{\stackup{\sim}{\operatorname{\longrightarrow}}}

\begin{document}


\title[Codescent theory II]{Codescent theory~II\,: Cofibrant approximations}

\author{Paul BALMER and Michel MATTHEY}

\address{Department of Mathematics, ETH Zentrum, CH-8092 Zürich, Switzerland}

\email{paul.balmer@math.ethz.ch}

\urladdr{http://www.math.ethz.ch/$\sim$balmer}

\address{Department of Mathematics, ETH Zentrum, CH-8092 Zürich, Switzerland}

\email{michel.matthey@math.ethz.ch}

\urladdr{http://www.math.ethz.ch/$\sim$matthey}

\thanks{Research supported by Swiss National Science Foundation, grant~620-66065.01}



\date{August 6, 2003}



\begin{abstract}
We establish a general method to produce cofibrant approximations in the model
category $\USCD{S}{C}{D}$ of $\cat{S}$-valued $\cat{C}$-indexed diagrams with $\cat{D}$-weak
equivalences and $\cat{D}$-fibrations. We also present explicit examples of such
approximations. Here, $\cat{S}$ is an arbitrary cofibrantly generated simplicial model
category and $\cat{D}\subset\cat{C}$ are small categories. An application to the notion of
homotopy colimit is presented.
\end{abstract}


\maketitle


\section{Introduction}

\label{s-IntroII}


%
The present paper may be read independently of Part I.

It is an important problem to understand model structures on categories of diagrams $\cat{S}^{\cat{C}}$,
indexed by a small category $\cat{C}$ and with values in some model category~$\cat{S}$, like for instance
the category of (compactly generated Hausdorff) topological spaces. One recent illustration of this importance,
among many others, is given in~\cite{bamaR}, where we show that the $K$-theoretic Isomorphism Conjectures
boil down to understanding cofibrant approximations in a suitable category of diagrams. In this spirit,
cofibrant approximations might be thought of as global assembly maps.

For an arbitrary model category $\cat{S}$, there is in general no known model structure on $\cat{C}$-indexed
diagrams $\cat{S}^{\cat{C}}$, with objectwise weak equivalences. Hence the notion of ``model approximation'' of Chach{\'o}lski
and Scherer~\cite{chsch}, that is not used here. Nevertheless, it is well-known that if the model category $\cat{S}$
is \emph{cofibrantly generated}, one can create a so-called \emph{left model structure} on $\cat{S}^{\cat{C}}$,
denoted $\USC{S}{C}$ hereafter, by defining the weak-equivalences and the fibrations $\cat{C}$-objectwise and
by forcing the cofibrations by a left lifting property. In Part~I, we even needed a \emph{$\cat{D}$-relative}
model structure on $\cat{S}^{\cat{C}}$, which we have denoted $\USCD{S}{C}{D}$, where the weak-equivalences and
the fibrations are defined $\cat{D}$-objectwise only, on a subset of objects $\cat{D}\subset\cat{C}$,
see~\cite[Thm~\ref{CDmodel-thm}]{bamaI}. The reader of Dugger~\cite{dugg} and Hirschhorn~\cite{hirsch} can
as well keep the absolute case $\cat{D}=\cat{C}$ in mind.

Unfortunately, these left model constructions, although very popular, usually leave the \emph{cofibrations} mysterious.
The factorization axiom in $\cat{S}^{\cat{C}}$, which guarantees their abundance, generally roots back to Quillen's small object
argument and this makes it hard to control what cofibrant approximations are. To be on the safe side, recall the terminology\,:
a weak equivalence $\xi\colon\mathcal{Q}X\too X$ with $\mathcal{Q}X$ cofibrant is called \emph{a cofibrant approximation}
of $X$.

In this paper, assuming that $\cat{S}$ is a cofibrantly generated \emph{simplicial} model category, we give an explicit
construction of cofibrant approximations in categories of $\cat{S}$-valued diagrams. We do this in the relative
case $\USCD{S}{C}{D}$ as well, mainly because we need it in~\cite{bamaR}.

\goodbreak

Of course, for $\cat{C}$ reduced to a point, $\cat{S}^{\cat{C}}$ is nothing but $\cat{S}$ and there is no hope that a general
process for $\cat{C}$-indexed diagrams suddenly provides us with new cofibrant approximations in an arbitrary $\cat{S}$.
Therefore, our method focusses on the ``diagrammatic part'' of the story and we consider cofibrant
approximations in $\cat{S}$ itself as being under control. This cofibrant approximation in the category $\cat{S}$
could even be the identity if everybody is cofibrant in $\cat{S}$, like \eg in the category $\sSets$ of simplicial sets.

Let $X\in\cat{S}^{\cat{C}}$ be a diagram. Let us start looking for a cofibrant approximation $\xi\colon Y\to X$. The first
observation is that a cofibrant diagram $Y\in\cat{S}^{\cat{C}}$ is always objectwise cofibrant, \ie $Y(c)$
is cofibrant in $\cat{S}$ for all $c\in\cat{C}$. So, it can not harm to first replace $X$ objectwise by a functorial cofibrant
approximation in $\cat{S}$. We produce in this way a first cheap approximation $\eta:qX\to X$, where
$qX(c)=\mathcal{Q}_{\cat{S}}(X(c))$. Although $qX$ is objectwise cofibrant and the map $\eta$ is a $\cat{C}$-weak
equivalence, this $qX$ is in general not cofibrant as a diagram\,! So far, we have done as much as we could do just
using the category of values $\cat{S}$. We now need to turn to the internal structure of the index-category $\cat{C}$ and,
in the relative case, of the subcategory $\cat{D}\subset\cat{C}$.

Recall that a simplicial model category (see Appendix~\ref{AppIIA}) is in particular equipped with an action
$\odot\colon\sSets\times\cat{S}\too\cat{S}$ of the category of simplicial sets. This can be jazzed up into a ``tensor product''
$$
-\,\udCDdot{C}{D}{\otimes}\,-\ \colon\quad\sSets^{\cat{D}\op\times\cat{C}}\,\times\,\cat{S}^{\cat{D}}\ \too\ \cat{S}^{\cat{C}}\,.
$$
This is probably well-known to the experts but we shall explain this carefully in~Section~\ref{s-Couplings}.

Our point is that finding cofibrant
approximations of $\cat{C}$-indexed diagrams with values in $\cat{S}$ amounts to finding \emph{one} cofibrant
approximation of \emph{one} very special and canonical diagram $\bbF$, living in $\sSets^{\cat{D}\op\times\cat{C}}$ and described
below, and then to tensor it with any object we want to cofibrantly approximate. In some sense, this diagram $\bbF$
encodes the purely $(\cat{C},\cat{D})$-part of the problem. We state
the following main result in relative form and then unfold the case $\cat{C}=\cat{D}$.

\begin{Thm}
\label{main-approx-thm}
Let $\cat{C}$ be a small category and let $\cat{D}\subset\cat{C}$ be a full subcategory. Define
$\bbF\in\sSets^{\cat{D}\op\times\cat{C}}$ as follows\,: for any $d\in\cat{D}$ and $c\in\cat{C}$,
$$
\bbF(d,c):=\mor_{\cat{C}}(d,c)\quad\in\,\sSets\,,
$$
where the set $\mor_{\cat{C}}(d,c)$ is seen as a constant (or discrete) simplicial set. Choose one
cofibrant approximation of $\,\bbF$ in $\USCD{\sSets}{D\op\times\cat{C}\,}{\,\cat{D}\op\times\cat{D}}$\,:
$$
\bbE\too\bbF\,.
$$
Let now $\cat{S}$ be a cofibrantly generated simplicial model category and let $X\in\cat{S}^{\cat{C}}$ be
a diagram. Let $qX\to X$ be an objectwise cofibrant approximation of $X$. Then
$$
\xi_{X}\,\colon\quad{\bbE}\;\;\udCDdot{C}{D}{\otimes}\;\res_{\cat{D}}^{\cat{C}}\,(qX)\ \too\ X
$$
is a cofibrant approximation of $X$ in the model category $\USCD{S}{C}{D}$. In particular, if $X\in\cat
{S}^{\cat{C}}$ is objectwise cofibrant, then $\,{\bbE}\,\udCDdot{D}{C}{\otimes}\res_{\cat{D}}^{\cat{C}}X\,
\too X$ is a cofibrant approximation of $X$ in $\USCD{S}{C}{D}$. The morphism $\xi_{X}$ is the ``obvious''
one (see Theorem~\ref{cof-approx-USCD-thm}).
\end{Thm}

We shall need in~\cite{bamaR} the above generality. However, since the absolute case $\cat{D}=\cat{C}$ widely
predominates in the literature, we now unfold the above theorem in this situation; moreover the morphism $\xi_{X}$
then becomes quite explicit.

\begin{Cor}
\label{main-approx-cor}
Let $\cat{C}$ be a small category. Consider $\bbF\in\sSets^{\cat{C}\op\times\cat{C}}$, defined for any
$c',c\in\cat{C}$ by
$
\bbF(c',c):=\mor_{\cat{C}}(c',c)\in\sSets
$
(simplicially constant). Choose a cofibrant approximation
$
\vartheta\colon\bbE\too\bbF
$
of $\bbF$ in the left model structure on $\sSets^{\cat{C}\op\times\cat{C}}$.
Let $\cat{S}$ be a cofibrantly generated simplicial model category and $X\in\cat{S}^{\cat{C}}$ a
diagram. Let $\eta\colon qX\to X$ be an objectwise cofibrant approximation of $X$. Then the composite
$$
\xi_{X}\colon
\xymatrix@C=1.3em{
{\bbE}\ \udCDdot{C}{C}{\otimes}\ qX
\ar[rrr]^-{\vartheta\,\udCDdot{C}{C}{\otimes}\,qX}
&&& {\bbF}\,\udCDdot{C}{C}{\otimes}qX
\cong qX \ar[rr]^-{\eta}
&& X}
$$
is a cofibrant approximation of $X$ in the left model category structure on $\cat{S}^{\cat{C}}$.
\end{Cor}

The customer is now entitled to ask for an explicit cofibrant approximation $\bbE\to \bbF$ of the object $\bbF\in
\sSets^{\cat{D}\op\times\cat{C}}$ of Theorem~\ref{main-approx-thm} with respect to the relative model structure
$\USCD{\sSets}{D\op\times\cat{C}\,}
{\,\cat{D}\op\times\cat{D}}$. Here it comes (see Theorem~\ref{EC-cofibrant-thm} below).

\begin{Thm}
\label{main-E-thm}
Let $\cat{C}$ be a small category and let $\cat{D}\subset\cat{C}$ be a full subcategory. Define $\bbF\in\sSets^{\cat
{D}\op\times\cat{C}}$ as above, \ie $\bbF(d,c)=\mor_{\cat{C}}(d,c)\in\sSets$ (simplicially constant). Define the
functor $\bbE\in\sSets^{\cat{D}\op\times\cat{C}}$ as follows. For $d\in\cat{D}$ and $c\in\cat{C}$ put
$$
\bbE(d,c):=B\,(d\comma\cat{D}\commaC{C}c)\op
$$
with $B$ standing for the usual nerve of the category $(d\comma\cat{D}\commaC{C}c)\op$; here $(d\comma\cat{D}\commaC{C}c)$
is the ``comma'' category of triples $(\alpha,x,\beta)$ where $x\in\cat{D}$ and where $\alpha\colon d\to x$ and
$\beta\colon x\to c$ are arrows in $\cat{C}$, the morphisms in $(d\comma\cat{D}\commaC{C}c)$ are the
obvious ones (see~\ref{commaC-not}). Consider the morphism $\vartheta
\colon\bbE\too \bbF$ given for $d\in\cat{D}$ and $c\in\cat{C}$ by the~map
$$
\vartheta(d,c)\colon\quad
\bbE(d,c)=B\,(d\comma\cat{D}\commaC{C}c)\op
\;\too\;
\mor_{\cat{C}}(d,c)=\bbF(d,c)\,,
$$
which is the evident ``composition of everything'' in each simplicial degree.
Then, $\vartheta\colon\bbE\too\bbF$ is a cofibrant approximation of $\bbF$ in $\USCD{\sSets}{D\op\times\cat{C}\,}
{\,\cat{D}\op\times\cat{D}}$.

\smallbreak

Similarly, there is a cofibrant approximation $\vthetapr\colon\deco{\bbE}\too\bbF$, where the object $\deco{\bbE}$ is defined
as above but without ``op'', that is, $\deco{\bbE}(d,c):=B\,(d\comma\cat{D}\commaC{C}c)$.
\end{Thm}

Combining Theorems~\ref{main-approx-thm} and~\ref{main-E-thm}, we obviously get explicit cofibrant approximations
in $\USCD{S}{C}{D}$, see~\ref{main-cor-def} below. Let us stress the universal character of these results.
Mastering one cofibrant approximation of one particular diagram $\bbF$ taking values in \emph{simplicial sets}
yields cofibrant approximations in $\cat{C}$-indexed diagrams with values in \emph{arbitrary reasonable model
categories}. Some authors, who like to think --~and to use~-- that whatever holds for diagrams in $\sSets$ remains
true for diagrams in familiar model categories, now have a rigorous argument at their disposal.

\smallbreak

As already observed in Hollender-Vogt~\cite{holvog} in the special case of topological spaces, the above diagram
$\bbE$ is of central importance for homotopy colimits.
The last short section of the paper is an application of the above to a question regarding homotopy
colimits, that
many topologists might have asked themselves once. Two approaches to homotopy colimits are available, both originating
from the work of Bousfield-Kan~\cite{boukan}. First, one can think of
the homotopy colimit basically as an esoteric but explicit formula which one can apply to whatever moves around, say, in
any category with an action of simplicial
sets. The second approach, slightly more conceptual, prefers to dwell on the problem that good old colimits \emph{do not} preserve
weak-equivalences, whereas homotopy colimits \emph{should ideally} preserve them. Homotopy colimits would therefore
be better suited for homotopy theory.
In a more model theoretical language, the
homotopy colimit should be thought of as the \emph{left derived functor}
of the colimit. Such a definition only makes sense if the category of diagrams is a model category, like for instance
if $\cat{S}$ is a cofibrantly generated model category. When both approaches make sense, namely if $\cat{S}$ is a
cofibrantly generated simplicial model category, as we consider here, it is reasonable and legitimate to ask whether both
notions coincide, \ie are weakly equivalent. The obvious
obstruction is of course that the standard homotopy colimit (the first mentioned above), if it ends up being
homotopy equivalent to a left derived functor, will itself preserve weak-equivalences. This is in fact
the one and only obstruction, as explained in Theorem~\ref{Lcolim=hocolim-thm} and Remark~\ref{zig-zag-rem},
which in particular imply the following result.

\begin{Thm}
\label{hocolim=Lcolim-thm}
Let $\cat{S}$ be a cofibrantly generated simplicial model category and let $\cat{C}$ be a small category. Assume
that $\hocolim_{\cat{C}}\colon \cat{S}^{\cat{C}}\too\cat{S}$ takes $\cat{C}$-weak equivalences
to weak-equivalences. Then there is a natural zig-zag of two weak-equivalences
$$
\hocolim_{\cat{C}}X\sim
\Lcolim_{\cat{C}}X
$$
for any $X\in\cat{S}^{\cat{C}}$, where $\Lcolim_{\cat{C}}X$ is any left derived functor of $\colim_{\cat{C}}$, that
is, $\Lcolim_{\cat{C}}=\colim_{\cat{C}}\circ \mathcal{Q}$ where $\mathcal{Q}$ is any functorial cofibrant approximation
in the category $\cat{S}^{\cat{C}}$ with the absolute model structure $\USC{S}{C}$.
Moreover, for $\mathcal{Q}$ suitably constructed with the methods of~\ref{main-approx-cor} and~\ref{main-E-thm}
above, this zig-zag of weak-equivalences reduces to a single natural weak-equivalence $\hocolim_{\cat{C}}X\simtoo
\Lcolim_{\cat{C}}X$.
\end{Thm}

\smallbreak

The organization of the article should be clear from this introduction and the table of contents below. Let us simply add
that we need some flexibility in ``couplings'' like the $\odot$ or the $\udCDdot{C}{D}{\otimes}$ considered above and
that this is better understood when abstracted a bit into a general notion of coupling $\cat{S}_{1}\times\cat{S}_{2}\too\cat
{S}_{3}$, with three possibly different categories involved. This gives us a fair chance to understand who does what in the
subsequent constructions. This formalism is developed in Sections~\ref{s-Couplings} and~\ref{s-couplings-model}.

Several special cases of our constructions are already available in the literature, like in the very complete~\cite{hirsch}
for instance.
Our main result which we use in the sequel~\cite{bamaR}, \ie the explicit description of cofibrant approximations
in the relative structure $\USCD{S}{C}{D}$, is clearly new. We also hope that the reader will benefit from the
systematic organization and from the relative concision of this article.


\tableofcontents
\goodbreak


\section{Notations, diagram-categories and recollection on $\USCD{S}{C}{D}$}

\label{s-Recoll-Part-I}


This short section introduces the notations used in the sequel and presents some
basic facts concerning categories of diagrams. It also contains a little summary
of what is needed from Part~I on the model category $\USCD{S}{C}{D}$.

\medskip

We refer to Mac Lane's book \cite{MacLane} for purely categorical questions.

\medskip

We refer to Hirschhorn~\cite{hirsch} or to Hovey~\cite{hov} for model category questions. Appendix~\ref{modcat-app} of
Part~I gives a concise list of prerequisites. The notion of simplicial model category being central here, it is recalled
in Appendix~\ref{AppIIA} of the present Part.

In a model category $\cat{M}$, recall the distinction between \emph{a} cofibrant \emph{approximation} of an object
$X\in\cat{M}$, meaning a weak-equivalence $\xi\colon\mathcal{Q}X\to X$ with $\mathcal{Q}X$ cofibrant, and \emph{the}
cofibrant \emph{replacement}, the latter being the cofibrant approximation obtained from the functorial factorization
axiom applied to the morphism
$\varnothing\to X$; this means in particular that $\xi$ is also a fibration in the latter case.

We call a cofibrant approximation $(\mathcal{Q},\xi)$ \emph{functorial}
if $\mathcal{Q}$ is a functor and if $\xi\colon\mathcal{Q}\too\id_{\cat{M}}$ is a natural transformation.
The cofibrant replacement is functorial and we sometimes designate it by $(Q_{\cat{M}},\xi^{\cat{M}})$.

\goodbreak \bigbreak
\centerline{*\ *\ *}
\smallbreak

\begin{Not}
\label{Cst-not}
Let $\cat{S}$ be a category and $\cat{C}$ a small category. Denote by $\cat{S}^{\cat{C}}$ the
category of functors from $\cat{C}$ to $\cat{S}$, also called \emph{$\cat{S}$-valued $\cat{C}$-indexed diagrams}. Write
$$
\Cst{\cat{S}}{\cat{C}}\colon\cat{S}\longrightarrow\cat{S}^{\cat{C}},
\quad
s\longmapsto\big(\cat{C}\rightarrow\cat{S},
\,c\mapsto s\big)
$$
for the ``constantification'' functor; this way, we can view $\cat{S}$ as a subcategory of $\cat{S}^
{\cat{C}}$. We also write $\underline{s}$ for $\Cst{\cat{S}}{\cat{C}}(s)$. Note that the colimit functor
$\colim_{\cat{C}}\colon\cat{S}^{\cat{C}}\longrightarrow\cat{S}$, when it exists, is left adjoint to the functor
$\Cst{\cat{S}}{\cat{C}}$.
\end{Not}

We denote by $\Sets$ the category of sets and by $\sSets=\Sets^{\DDelta\!\op}$ the category of
simplicial sets.

\begin{Rem}
\label{Cst-level-zero-rem}
The usual constantification functor $\Cst{\Sets}{\DDelta\!\op}\colon\Sets\longrightarrow\sSets$
has the following ``level-zero functor'' as right adjoint\,:
$$
(-)_{0}\colon\sSets\longrightarrow\Sets,\quad K=K_{\bullet}\longmapsto K_{0}\,.
$$
\end{Rem}

\begin{Conv}
We denote a functor $F$ also as $F(-)$ or $F(?)$. If in some ``formula'' two functors
$F$ and $G$ with the same source category $\cat{A}$ are involved, we write $F(?)$ and $G(?)$ to
stress the fact that we evaluate $F$ and $G$ at the same dummy-variable object $?$ of $\cat{A}$.
We adopt similar notations with $\qquest$ and $\qqquest$ in place of $?$, usually when several
dummy-variables are involved.
\end{Conv}

\begin{Not}
\label{switch-not}
For categories $\cat{S}$, $\cat{A}$ and $\cat{C}$, with $\cat{A}$ and $\cat{C}$ small,
we make the following obvious identifications of categories of diagrams\,:
$$
\big(\cat{S}^{\cat{C}}\big)^{\cat{A}}=\cat{S}^{\cat{C}\times\cat{A}}\qquad\mbox{and}\qquad
\big(\cat{S}^{\cat{A}}\big)^{\cat{C}}=\cat{S}^{\cat{A}\times\cat{C}}\,.
$$
We denote the evident ``switch functor'' by\ \
$
\sigma_{\cat{C},\cat{A}}\colon\cat{S}^{\cat{C}\times\cat{A}}\stackrel{\cong}{\longrightarrow}
\cat{S}^{\cat{A}\times\cat{C}}
$.
\end{Not}

\medbreak
\centerline{*\ *\ *}
\smallbreak

\begin{Def}
Let $\cat{S}$ be a cofibrantly generated model category, and $\cat{D}$ a subcategory of a small
category $\cat{C}$. A morphism $\varphi\colon X\longrightarrow Y$ in $\cat{S}^{\cat{C}}$ is
a \emph{$\cat{D}$-isomorphism} (resp.\ a \emph{$\cat{D}$-weak equivalence}, a \emph{$\cat{D}$-fibration}
or a \emph{trivial $\cat{D}$-fibration}) if for every object $d\in\cat{D}$, the morphism $\varphi(d)\colon X(d)
\longrightarrow Y(d)$ is an isomorphism (resp.\ a weak equivalence, a fibration or a trivial fibration)
in $\cat{S}$.
\end{Def}

\begin{PropDef}
\label{USCD-propdef}
Let $\cat{S}$ be a cofibrantly generated model category, and $\cat{D}$ a subcategory of a small
category $\cat{C}$. Then, there is a model category structure on the category $\cat{S}^{\cat{C}}$
of $\cat{C}$-indexed $\cat{S}$-valued diagrams, where the weak equivalences are the $\cat{D}$-weak equivalences
and the fibrations are the $\cat{D}$-fibrations. It is denoted by $\USCD{S}{C}{D}$ and is called
the \emph{$\cat{D}$-relative model structure} on $\cat{S}^{\cat{C}}$. If a diagram $X\in\cat{S}^{\cat{C}}$
is cofibrant in $\USCD{S}{C}{D}$, we call it \emph{$\cat{D}$-cofibrant}. When $\cat{D}=\cat{C}$, we
also write $\USC{S}{C}$ for $\USCD{S}{C}{C}$ and call it the \emph{absolute model structure} on $\cat{S}^{\cat{C}}$.
\end{PropDef}

See \cite[Thm.\ \ref{CDmodel-thm}]{bamaI}. Note that the model category $\USCD{S}{C}{D}$ does not depend on the category
structure of $\cat{D}$,
but only on the underlying set of objects $\obj(\cat{D})$.

\begin{Rem}
\label{objwise-cof-rem}
Let $\cat{C}$ be a small category, $\cat{D}\subset\cat{C}$ a subcategory and $\cat{S}$ a cofibrantly generated model category.
It is proven in \cite[Prop.~\ref{Dcof-casc-prop}]{bamaI} (see also Remark~\ref{cof-objw-rem}
therein) that if $f\colon X_{1}\longrightarrow X_{2}$
is a cofibration in $\USCD{S}{C}{D}$, then it is $\cat{C}$-objectwise a cofibration, that is,
$f(c)\colon X_{1}(c)\longrightarrow X_{2}(c)$ is a cofibration in $\cat{S}$ for each $c\in
\cat{C}$. In particular, if $X\in\USCD{S}{C}{D}$ is a cofibrant diagram, then $X(c)\in\cat{S}$ is cofibrant
as well, for every object $c\in\cat{C}$ (and not just $\cat{D}$-objectwise).
\end{Rem}

Recall that a category is \emph{(co)complete} if it admits all small (co)limits.

\begin{Rem}
\label{res-ind-Q-adj-rem}
Let $\cat{S}$ be a cocomplete category and $\cat{D}$ a subcategory of a small category $\cat{C}$.
The restriction functor $\res_{\cat{D}}^{\cat{C}}$, from $\cat{S}^{\cat{C}}$ to $\cat{S}^{\cat{D}}$, has a
left adjoint
$$
\ind_{\cat{D}}^{\cat{C}}\colon\cat{S}^{\cat{D}}\adjtoo
\cat{S}^{\cat{C}}\noloc\res_{\cat{D}}^{\cat{C}}
$$
called the induction functor (see for instance~\cite[App.~\ref{Kan-app}]{bamaI}). Assume now that $\cat{S}$ is a
cofibrantly generated
model category. It is clear that the functor $\res_{\cat{D}}^{\cat{C}}$ takes (trivial) $\cat{D}$-fibrations
to (trivial) $\cat{D}$-fibrations, so that we have a Quillen adjunction
$$
\ind_{\cat{D}}^{\cat{C}}\colon\USC{S}{D}\adjtoo\USCD{S}{C}{D}\noloc\res_{\cat{D}}^{\cat{C}}\,.
$$
\end{Rem}

The forthcoming two observations will be repeatedly used in the sequel.

\begin{Rem}
\label{playing-U-U-rem}
For a cofibrantly generated model category $\cat{S}$ and two small categories $\cat{A}$
and $\cat{C}$, the identification $\big(\cat{S}^{\cat{C}}\big)^{\cat{A}}=\cat{S}^{\cat{C}
\times\cat{A}}$ of Notation \ref{switch-not} is not just an identification of mere categories,
but really an identification of model categories (up to the choice of the functorial
factorizations)\,:
$$
\UUSC{\USC{S}{C}}{\cat{A}}=\UUSC{\cat{S}}{\cat{C}\times\cat{A}}
\qquad\mbox{or even}\qquad
\USCD{{}\USCD{S}{C}{D}}{A}{B}=\USCD{S}{C\times\cat{A}}{D\times\cat{B}}
$$
for any subcategories $\cat{D}\subset\cat{C}$ and $\cat{B}\subset\cat{A}$ (the latter $\cat{B}\subset\cat{A}$ is not used below).
\end{Rem}

In the next remark, we use the opposite category $\cat{A\op}$ of $\cat{A}$, instead of $\cat{A}$ directly,
only for cosmetic reasons justified by the use we make of the remark later on.

\begin{Rem}
\label{CD-obj-wise-cof-rem}
Let $\cat{A}$ and $\cat{D}\subset\cat{C}$ be small categories and let $\cat{S}$ be a cofibrantly generated
model category. Suppose that $\mathcal{X}(\qquest,?)$ is a cofibrant diagram in the model
category $\USCD{S}{A\op\times\cat{C}}{A\op\times\cat{D}}$. Then, $\mathcal{X}(\qquest,?)$ is $\cat{A}\op$-objectwise,
$\cat{C}$-objectwise and $\cat{A}\op\times\cat{C}$-objectwise cofibrant. More explicitly, for objects
$a\in\cat{A}$ and $c\in\cat{C}$, the three objects $\mathcal{X}(a,?)\in\USCD{S}{C}{D}$, $\mathcal{X}
(\qquest,c)\in\USC{S}{A\op}$ and $\mathcal{X}(a,c)\in\cat{S}$ are cofibrant, as follows from
Remark~\ref{objwise-cof-rem} (cf.\ \ref{switch-not} and \ref{playing-U-U-rem}).
\end{Rem}


\section{Couplings of categories and couplings of diagrams}

\label{s-Couplings}


In this section, we consider functors $\odot\colon\cat{S}_{1}\times\cat{S}_{2}\longrightarrow
\cat{S}_{3}$ from a product of two categories to a possibly different category; we call them
couplings. This is in particular studied when some of the categories $\cat{S}_{i}$ are replaced by the
category of $\cat{S}_{i}$-valued functors on a given small category. More explicitly,
letting $\cat{C}$ and $\cat{A}$ be small categories, we will induce up two couplings
$$
\uCdot{\pair{C}{\;\,}}{\odot}\colon
\big(\cat{S}_{1}\big)^{\cat{C}}\times\cat{S}_{2}
\longrightarrow\big(\cat{S}_{3}\big)^{\cat{C}}
\qquad\mbox{and}\qquad
\dCdot{A}{\otimes}\colon
\big(\cat{S}_{1}\big)^{\cat{A}\op}\times\big(\cat{S}_{2}\big)^{\cat{A}}
\longrightarrow\cat{S}_{3}\,.
$$
Combining both constructions, we will deduce yet another useful coupling
$$
\udCDdot{C}{A}{\otimes}\colon\big(\cat{S}_{1}\big)^{\cat{A}\op\times\cat{C}}\times\big(\cat
{S}_{2}\big)^{\cat{A}}\longrightarrow\big(\cat{S}_{3}\big)^{\cat{C}}\,.
$$

The following adjunctions are our basic tools.

\begin{Def}
\label{coupling-def}
Let $\cat{S}_{1}$, $\cat{S}_{2}$ and $\cat{S}_{3}$ be three categories.
\begin{itemize}
\item [(i)] An \emph{$\cat{S}_{3}$-valued coupling of $\cat{S}_{1}$ with $\cat{S}_{2}$}
is a (bi-)functor
$$
\itemspace\odot\colon\cat{S}_{1}\times\cat{S}_{2}\longrightarrow\cat{S}_{3},\quad(x,y)\longmapsto
x\odot y\,.
$$
\item [(ii)] The coupling $\odot$ is called \emph{right tensorial} if there exists a functor
$$
\itemspace\rmap{\odot}\colon\cat{S}_{2}\op\times\cat{S}_{3}\longrightarrow\cat{S}_{1},\quad(y,z)\longmapsto
\rmap{\odot}(y,z)
$$
such that for every $x\in\cat{S}_{1}$, $y\in\cat{S}_{2}$ and $z\in\cat{S}_{3}$, there
is a natural bijection
$$
\itemspace\mor_{\cat{S}_{3}}(x\odot y,z)\cong\mor_{\cat{S}_{1}}\!\big(x,\rmap{\odot}(y,z)\big)\,,
$$
in other words, if for every $y\in\cat{S}_{2}$, there is an adjunction
$$
\itemspace(-)\odot y\colon\cat{S}_{1}\adjtoo\cat{S}_{3}\noloc\rmap{\odot}(y,-)
$$
depending functorially on $y$; we call $\rmap{\odot}$ the \emph{right mapping functor} of~$\odot$.
\itemspace
\item [(iii)] The coupling $\odot$ is called \emph{left tensorial} if there exists a functor
$$
\itemspace\lmap{\odot}\colon\cat{S}_{1}\op\times\cat{S}_{3}\longrightarrow\cat{S}_{2},\quad(x,z)\longmapsto
\lmap{\odot}(x,z)
$$
such that for every $x\in\cat{S}_{1}$, $y\in\cat{S}_{2}$ and $z\in\cat{S}_{3}$, there
is a natural bijection
$$
\itemspace\mor_{\cat{S}_{3}}(x\odot y,z)\cong\mor_{\cat{S}_{2}}\!\big(y,\lmap{\odot}(x,z)\big)\,,
$$
in other words, if for every $x\in\cat{S}_{1}$, there is an adjunction
$$
\itemspace x\odot(-)\colon\cat{S}_{2}\adjtoo\cat{S}_{3}\noloc\lmap{\odot}(x,-)
$$
depending functorially on $x$; we call $\lmap{\odot}$ the \emph{left mapping functor} of
$\odot$.
\end{itemize}
\end{Def}

Of course, for a right (resp.\ left) tensorial coupling the right (resp.\ left) mapping
functor is uniquely determined up to unique natural isomorphism.

\begin{Ex}
\label{coupling-sets-ex}
If $\cat{S}$ is a category with small coproducts (for example a cocomplete category),
then there is a canonical $\cat{S}$-valued coupling of $\Sets$ with $\cat{S}$ given by
$$
\boxdot\colon\;\Sets\times\cat{S}\longrightarrow\cat{S},\qquad(K,s)\longmapsto
K\boxdot s:=\colim_{K}\,\underline{s}
=\coprod_{K}s\,,
$$
where in the indicated colimit, we view $K$ as a discrete category, that is, with $K$ as
set of objects and only with identity morphisms, and $\underline{s}$ denotes $\Cst{\cat{S}}
{K}(s)$. For $K\in\Sets$ and $s,s'\in\cat{S}$, the universal property of the coproduct provides
a natural bijection
$$
\mor_{\cat{S}}\!\big({\textstyle \coprod_{K}s},s'\big)\cong\mor_{\Sets}\!\big(K,\mor_{\cat{S}}(s,s')\big)\,.
$$
This shows that the coupling $\boxdot$ is right tensorial with right mapping functor $\mor_{\cat{S}}$.
Note that for the point $*\in\Sets$ and any object $s\in\cat{S}$, we have a canonical and natural
isomorphism $*\boxdot s\cong s$.
\end{Ex}

\begin{Ex}
\label{coupling-ex}
If $\cat{S}$ is a simplicial model category (see \ref{simpl-mod-cat-def}), then the ``action''
$\odot\colon\sSets\times\cat{S}\too\cat{S}$ (see \ref{odot-nat}) of the category $\sSets$ of simplicial
sets on $\cat{S}$ is an $\cat{S}$-valued coupling of $\sSets$ with $\cat{S}$. It is right tensorial with
$\rmap{\odot}(y,z):=\Map(y,z)$ for every $y,z\in\cat{S}$, by virtue of Axiom~\ref{simpl-mod-cat-def}\,(2),
using also \ref{simpl-mod-cat-def}\,(iii). It is left tensorial with $\lmap{\odot}(x,z)=z^{x}$ for every
$x\in\sSets$ and $z\in\cat{S}$, by virtue of Axiom~\ref{simpl-mod-cat-def}\,(3).
\end{Ex}

It turns out that the couplings of Examples \ref{coupling-sets-ex} and \ref{coupling-ex}
are ``compatible''.

\begin{Lem}
\label{compatibility-lem}
Let $\cat{S}$ be a simplicial model category. Consider the couplings $\boxdot$ of
Example~\ref{coupling-sets-ex} and $\odot$ of Example~\ref{coupling-ex}. Then,
the following diagram commutes up to natural isomorphism
$$
\xymatrix @R=1.5em{
{\Sets\times\cat{S}} \ar[rr]^-{\Cst{\Sets}{\DDelta\!\op}\times\cat{S}} \ar[rd]_{\boxdot} & &
{\sSets\times\cat{S}} \ar[ld]^{\odot} \\
& \cat{S} & \\
}
$$
In particular, for the point $*=\Delta^{0}\in\sSets$ and an arbitrary object $s\in\cat{S}$,
there is a natural isomorphism $*\odot s\cong s$.
\end{Lem}

\Prf
By right adjunction, it suffices to prove for all $s\in\cat{S}$ the commutativity of
$$
\xymatrix @R=1.5em{
{\Sets} & & {\sSets} \ar[ll]_-{(-)_{0}} \\
& \cat{S} \ar[lu]^{\mor_{\cat{S}}(s,-)=\rmap{\,\boxdot}(s,-)\quad} \ar[ru]_{\quad\rmap{\odot}(s,-)=\Map(s,-)} & \\
}
$$
which is precisely one of the axioms of simplicial model categories\,: see \ref{simpl-mod-cat-def}\,(iii).
\qed

\begin{Rem}
If $\odot\colon\cat{S}_{1}\times\cat{S}_{2}\too\cat{S}_{3}$ is a coupling, then the composite functor
$$
\xymatrix @C=3.5em{
{}^{t}\!\odot\colon\cat{S}_{2}\times\cat{S}_{1} \ar[r]^-{{\rm switch}} & \cat{S}_{1}\times\cat{S}_{2}
\ar[r]^-{\odot} & \cat{S}_{3},\quad y\,{}^{t}\!\!\odot x:=x\odot y
}
$$
is also a coupling. It is left (resp.\ right) tensorial if and only if $\odot$ is right (resp.\ left) tensorial, with
$\lmap{{}^{t}\!\nnspace\odot}(?,\qquest):=\rmap{\odot}(\qquest,?)$ (and similarly in the other case). Using this
observation, we shall focus on right tensorial couplings in the sequel and leave the dual statements to the reader.
\end{Rem}

\medbreak
\centerline{*\ *\ *}
\medbreak

We pass to the first construction of induced couplings on categories of diagrams.

\begin{Lem}
\label{odot-induced1-lem}
Let $\odot\colon\cat{S}_{1}\times\cat{S}_{2}\too\cat{S}_{3}$ be a coupling and
let $\cat{C}$ be a small category. Then, the assignment
$
\uCdot{\pair{C}{C}}{\odot}\colon\cat{S}_{1}^{\cat{C}}\times\cat{S}_{2}^{\cat{C}}\longrightarrow
\cat{S}_{3}^{\cat{C}},\quad X\uCdot{\pair{C}{C}}{\odot}Y\;(?):=X(?)\odot Y(?)
$
is a coupling and so is its restriction (called the \emph{induced coupling over $\cat{C}$})
$$
\uCdot{\pair{C}{\;\,}}{\odot}\colon\cat{S}_{1}^{\cat{C}}\times\cat{S}_{2}\longrightarrow\cat{S}_{3}^{\cat{C}},
\qquad X\uCdot{\pair{C}{\;\,}}{\odot}y\;(?):=X(?)\odot y\,.
$$
If moreover $\odot$ is right tensorial, then $\uCdot{\pair{C}{\;\,}}{\odot}$ is a right tensorial coupling
with right mapping functor given by
$$
\rmap{\smalluCdot{\pair{C}{\;\,}}{\odot}}\colon\cat{S}_{2}\op\times\cat{S}_{3}^{\cat{C}}\longrightarrow
\cat{S}_{1}^{\cat{C}},
\qquad\rmap{\smalluCdot{\pair{C}{\;\,}}{\odot}}(y,Z)\,(?):=\rmap{\odot}\!\big(y,Z(?)
\big)\,.
$$
\end{Lem}

\Prf
It is plain that $\uCdot{\pair{C}{C}}{\odot}$ and $\uCdot{\pair{C}{\;\,}}{\odot}$ are couplings and let us prove the moreover
part. Let $X\in\cat{S}_{1}^{\cat{C}}$, $y\in\cat{S}_{2}$,
$Z\in\cat{S}_{3}^{\cat{C}}$.
By adjunction~\ref{coupling-def}\,(ii), for a given morphism $\alpha\colon a\longrightarrow b$ in $\cat{C}$, commutative
diagrams in $\cat{S}_{3}$ like $(*)$
$$
(*)\qquad
\vcenter{\xymatrix@R=1.7em{
X(a)\odot y \ar[r]^-{\psi_{a}} \ar[d]_{X(\alpha)\odot y} & Z(a) \ar[d]^{Z(\alpha)} \\
X(b)\odot y \ar[r]^-{\psi_{b}} & Z(b) \\
}}
\qquad
\vcenter{\xymatrix@R=1.7em{
X(a) \ar[r]^-{\tau_{a}} \ar[d]_{X(\alpha)} & \rmap{\odot}(y,Z(a)) \ar[d]^{\rmap{\odot}(y,Z(\alpha))} \\
X(b) \ar[r]^-{\tau_{b}} & \rmap{\odot}(y,Z(b)) \\
}}
\qquad(**)
$$
are in one-one correspondence with commutative diagrams in $\cat{S}_{1}$ like $(**)$. The set
$$
\left.\left\{(\psi_{c})_{c\in\cat{C}}\in\prod_{c\in\cat{C}}\mor_{\cat{S}_{3}}\!\big(X(c)\odot y,Z(c)\big)\,
\right|\,\mbox{$(*)$ commutes $\forall\alpha\colon a\rightarrow b$ in $\cat{C}$}\right\}\,,
$$
which coincides with the set $\mor_{\cat{S}_{3}^{\cat{C}}}\!\big(X\ \uCdot{\pair{C}{\;\,}}{\odot}y,Z\big)$, is
thus in bijection with the set
$$
\left.\left\{(\tau_{c})_{c\in\cat{C}}\in\prod_{c\in\cat{C}}\mor_{\cat{S}_{3}}\!\Big(X(c),\rmap{\odot}
\!\big(y,Z(c)\big)\Big)\,\right|\,\mbox{$(**)$ commutes $\forall\alpha\colon a\rightarrow b$ in $\cat{C}$}
\right\}\,,
$$
which is equal to $\mor_{\cat{S}_{1}^{\cat{C}}}\!\big(X,\rmap{\smalluCdot{\pair{C}{\;\,}}{\odot}}\!(y,Z)\big)$. The rest is routine.
\qed

\begin{Rem}
Note that in Lemma \ref{odot-induced1-lem}, when $\odot$ is right tensorial, it is generally
not true that $\uCdot{\pair{C}{C}}{\odot}$ is right tensorial as well, with a right mapping functor
which would be given, for an object $c\in\cat{C}$, by $\rmap{\smalluCdot{\pair{C}{C}}{\odot}}(Y,Z)\,(c)
:=\rmap{\odot}\!\big(Y(c),Z(c)\big)$; indeed, there is no reasonable way of making a functor
out of this objectwise assignment $\rmap{\smalluCdot{\pair{C}{C}}{\odot}}$ (unless $\cat{C}$
is discrete or $\cat{S}_{2}$ is trivial).
\end{Rem}

\begin{Not}
\label{other-induced-not}
Let $\odot\colon\cat{S}_{1}\times\cat{S}_{2}\too\cat{S}_{3}$ be a coupling and
$\cat{C}$ a small category. The coupling $\uCdot{\pair{C}{C}}{\odot}$ of Lemma \ref{odot-induced1-lem} has a second
 possible restriction, denoted by
$$
\uCdot{\pair{\;\,}{C}}{\odot}\colon\cat{S}_{1}\times\cat{S}_{2}^{\cat{C}}\longrightarrow\cat{S}_{3}^{\cat{C}},
\qquad x\uCdot{\pair{\;\,}{C}}{\odot}Y\,(?):=x\odot Y(?)\,.
$$
\end{Not}

\medbreak
\centerline{*\ *\ *}
\medbreak

Next, we study a second way to induce up couplings on categories of diagrams.

Recall from \cite[pp.\ 64--65]{MacLane} the notion of co-equalizer of a pair
$\xymatrix @C=1.5em{s_{1} \ar@<-.22em>[r]_-{g} \ar@<.22em>[r]^-{f} & s_{2}}$
of morphisms in a category $\cat{S}$, with same source and same target\,: it
is defined by
$$
\coeq\Big(\xymatrix @C=1.5em{s_{1} \ar@<-.22em>[r]_-{g} \ar@<.22em>[r]^-{f} & s_{2}}\Big)
:=\colim_{\coeqcat}\Big(\xymatrix @C=1.5em{s_{1} \ar@<-.22em>[r]_-{g} \ar@<.22em>[r]^-{f} &
s_{2}}\Big)\,,
$$
the universal morphism (or only its target) out of $s_{2}$ which ``co-equalizes'' $f$ and $g$.
\goodbreak

\begin{Def}
\label{tensor-def}
Let $\odot\colon\cat{S}_{1}\times\cat{S}_{2}\too\cat{S}_{3}$ be a coupling, with
$\cat{S}_{3}$ cocomplete, and let $\cat{A}$ be a small category. The \emph{tensor product
over $\cat{A}$ associated to $\odot$} is the assignment
$$
\dCdot{A}{\otimes}\colon\cat{S}_{1}^{\cat{A}\op}\times\cat{S}_{2}^{\cat{A}}\longrightarrow\cat{S}_{3},
\qquad(X,Y)\longmapsto X\dCdot{A}{\otimes}Y
$$
defined for $X\in\cat{S}_{1}^{\cat{A}\op}$ and $Y\in\cat{S}_{2}^{\cat{A}}$ by the formula
$$
X\dCdot{A}{\otimes}Y
=\coeq\Bigg(\;\coprod_{(b\stackrel{\!\alpha}{\rightarrow}b')\in
\cat{A}}X(b')\odot Y(b)
\;\;{\renewcommand{\arraystretch}{0.5}
\begin{array}{@{}c@{}}
\longrightarrow \\
\longrightarrow \\
\end{array}
\renewcommand{\arraystretch}{1}}\;\;
\coprod_{a\in\cat{A}}X(a)\odot Y(a)\;\Bigg)
$$
where the two indicated morphisms inside the co-equalizer are given, on a summand indexed by $\alpha\colon b\longrightarrow b'$,
by the compositions
$$
\xymatrix @C=2em @R=.1em {
&&& X(b)\odot Y(b) \ar[r]
& {\displaystyle \coprod_{a\in\cat{A}}^{\vphantom{\cat{A}}}}X(a)\odot Y(a)
\\
X(b')\odot Y(b)
\ar[rrru]^-{X(\alpha)\odot Y(b)\ \ }
\ar[rrrd]_-{X(b')\odot Y(\alpha)\ \ }
\\
&&& X(b')\odot Y(b') \ar[r]
& {\displaystyle \coprod_{a\in\cat{A}}^{\vphantom{\cat{A}}}}X(a)\odot Y(a)}
$$
It is sometimes useful to write $X\dCdot{A}{\otimes}Y$ as $X(?)\dCdot{A}{\otimes}Y(?)$ or even
as $X(?)\dCdot{?\in A}{\otimes}Y(?)$.
\end{Def}

\begin{Ex}
\label{colim=*tensor(-)-ex}
Let $\odot\colon\cat{S}_{1}\times\cat{S}_{2}\longrightarrow\cat{S}_{2}$ be a coupling (note that $\cat{S}_{3}=\cat{S}_{2}$).
That is,
$\odot$ is an ``action'' of a category $\cat{S}_{1}$ on $\cat{S}_{2}$. Suppose that $\cat{S}_{2}$ is
cocomplete. Assume that there exists an object $x_{0}\in\cat{S}_{1}$ such that the ``action''
of $x_{0}$ is trivial, more precisely, such that the functor $x_{0}\odot(-)$ is
isomorphic to the identity functor of $\cat{S}_{2}$. Then, for every small category $\cat{A}$,
there is an isomorphism of functors
$$
\underline{x_{0}}\,\dCdot{A}{\otimes}(-)\cong\colim_{\cat{A}}(-)\colon
\quad\cat{S}_{2}^{\cat{A}}\longrightarrow
\cat{S}_{2}\,,
$$
where $\underline{x_{0}}\in\cat{S}_{1}^{\cat{A}\op}$ is the constant diagram with value $x_{0}\in\cat{S}_{1}$.
Indeed, for $Y\in\cat{S}_{2}^{\cat{A}}$, one obtains a natural isomorphism
$$
\underline{x_{0}}\dCdot{A}{\otimes}Y\cong\colim_{\coeqcat}\left(\;\coprod_{(b\stackrel{\!\alpha}{\rightarrow}b')\in
\cat{A}}Y(b)
\xymatrix @C=3.8em{
{} \ar@<-.22em>[r]_-{\id} \ar@<.22em>[r]^-{\amalg_{\alpha}Y(\alpha)} & {}
}
\coprod_{a\in\cat{A}}Y(a)\right)\,,
$$
and the latter is easily seen to satisfy the universal property of the colimit of $Y$ over $\cat{A}$.
A typical situation where this applies is for $\odot$ as in Example~\ref{coupling-ex}, with $\cat{S}_{1}
:=\sSets$, $\cat{S}_{2}:=\cat{S}$ (a simplicial model category) and with $x_{0}:=\Delta^{0}\in\sSets$.
\end{Ex}

\begin{Lem}
\label{odot-induced2-lem}
Let $\odot\colon\cat{S}_{1}\times\cat{S}_{2}\too\cat{S}_{3}$ be a coupling, with
$\cat{S}_{3}$ cocomplete, and let $\cat{A}$ be a small category. Suppose that $\odot$ is
right tensorial. Then, the assignment
$$
\dCdot{A}{\otimes}\colon\cat{S}_{1}^{\cat{A}\op}\times\cat{S}_{2}^{\cat{A}}\longrightarrow\cat{S}_{3}
$$
is a right tensorial coupling with right mapping functor
$$
\rmap{\smalldCdot{A}{\otimes}}\colon
\big(\cat{S}_{2}^{\cat{A}}\big)\nnspace\op\times\cat{S}_{3}
\longrightarrow
\cat{S}_{1}^{\cat{A}\op},
\qquad
\rmap{\smalldCdot{A}{\otimes}}(Y,z)\,(?):=\rmap{\odot}\!\big(Y(?),z\big)\,.
$$
\end{Lem}

\Prf
The fact that $\dCdot{A}{\otimes}$ is a functor is a straightforward checking. Let $X\in\cat{S}_{1}^{\cat{A}\op}$, $Y\in
\cat{S}_{2}^{\cat{A}}$ and $z\in\cat{S}_{3}$.
By definition of co-equalizers, the set $\mor_{\cat{S}_{3}}(X\dCdot{A}{\otimes}Y,z)$ is equal to the set of those elements
$(\psi_{a})_{a\in\cat{A}}$ in $\prod_{a\in\cat{A}}\mor_{\cat{S}_{3}}\!\big(X(a)\odot
Y(a),z\big)$ such that for every morphism $\alpha\colon b\longrightarrow b'$ in $\cat{A}$,
the following left-hand (or equivalently right-hand) diagram commutes\,:
$$
\vcenter{\xymatrix @C=5em{
X(b')\odot Y(b) \ar[r]^-{X(b')\odot Y(\alpha)} \ar[d]^{X(\alpha)\odot Y(b)} & X(b')\odot Y(b')
\ar[d]^{\psi_{b'}} \\
X(b)\odot Y(b) \ar[r]^-{\psi_{b}} & z \\
}}
\qquad\mbox{or}\qquad
\vcenter{\xymatrix @C=4em @R=1.5em{
X(b')\odot Y(b') \ar[r]^-{\psi_{b'}} & z \ar@{=}[d] \\
X(b')\odot Y(b) \ar@{-->}[r] \ar[u]_{X(b')\odot Y(\alpha)} \ar[d]^{X(\alpha)\odot Y(b)} &
z \ar@{=}[d] \\
X(b)\odot Y(b) \ar[r]^-{\psi_{b}} & z \\
}}
$$
By adjunction~\ref{coupling-def}\,(ii), this set is in bijection with the set of those tuples $(\tau_{a})_{a\in\cat{A}}$
in $\prod_{a\in\cat{A}}
\mor_{\cat{S}_{1}}\!\big(X(a),\rmap{\odot}(Y(a),z)\big)$ making the following left-hand (or equivalently
right-hand) diagram commutative for every morphism $\alpha\colon b\longrightarrow b'$ in $\cat{A}$\,:
$$
\vcenter{\xymatrix @C=4em @R=1.5em{
X(b') \ar[r]^-{\tau_{b'}} \ar@{=}[d] & \rmap{\odot}(Y(b'),z) \ar[d]_{\rmap{\odot}(Y(\alpha),z)\!} \\
X(b') \ar@{-->}[r] \ar[d]_{X(\alpha)} &
\rmap{\odot}(Y(b),z) \ar@{=}[d] \\
X(b) \ar[r]^-{\tau_{b}} & \rmap{\odot}(Y(b),z) \\
}}
\qquad\mbox{or}\qquad
\vcenter{\xymatrix{
X(b') \ar[r]^-{\tau_{b'}} \ar[d]_{X(\alpha)} & \rmap{\odot}(Y(b'),z) \ar[d]_{\rmap{\odot}(Y(\alpha),z)\!} \\
X(b) \ar[r]^-{\tau_{b}} & \rmap{\odot}(Y(b),z) \\
}}
$$
This set is nothing but $\mor_{\cat{S}_{1}^{\cat{A}\op}}\!\big(X(?),\rmap{\odot}(Y(?),z)\big)$, as
was to be shown, the required naturality properties being, again, routine.
\qed

\goodbreak \bigbreak
\centerline{*\ *\ *}
\bigbreak

Mixing the two ways of inducing couplings seen so far, we have the immediate\,:

\medbreak
\begin{PropDef}
\label{bi-tensor-propdef}
Let $\odot\colon\cat{S}_{1}\times\cat{S}_{2}\too\cat{S}_{3}$ be a coupling, with
$\cat{S}_{3}$ cocomplete. Consider two small categories $\cat{A}$ and $\cat{C}$. Applying successively
the constructions of Lemma~\ref{odot-induced2-lem} and of Lemma~\ref{odot-induced1-lem}, we obtain a coupling
$$
\uCdot{\pair{C}{\quad}}{\raisebox{-.25em}{\big(}\nnspace\dCdot{A}{\otimes}\!\raisebox{-.25em}{\big)}}
\colon\cat{S}_{1}^{\cat{A}\op\times\cat{C}}\times\cat{S}_{2}^{\cat{A}}\longrightarrow\cat{S}_{3}^{\cat{C}}
$$
that we simply denote by $\udCDdot{C}{A}{\otimes}$.
If moreover the coupling $\odot$ is right tensorial, then so is $\udCDdot{C}{A}{\otimes}$ and its right mapping functor reads
$$
\rmap{\smalludCDdot{C}{A}{\otimes}}\colon
\big(\cat{S}_{2}^{\cat{A}}\big)\nnspace\op\times\cat{S}_{3}^{\cat{C}}
\longrightarrow\cat{S}_{1}^{\cat{A}\op\times\cat{C}},
\qquad\rmap{\smalludCDdot{C}{A}{\otimes}}(Y,Z)\,(\qquest,?):=
\rmap{\odot}\!\big(Y(\qquest),Z(?)\big)\,.
$$
We call $\udCDdot{C}{A}{\otimes}$ the \emph{bi-tensor product over $(\cat{C},\cat{A})$ induced by $\odot$}.
\qed
\end{PropDef}

Explicitly, for
$\mathcal{X}\in\cat{S}_{1}^{\cat{A}\op\times\cat{C}}$, for $Y\in\cat{S}_{2}^{\cat{A}}$ and for $c\in\cat{C}$, we have
$$
\arraycolsep1pt
\renewcommand{\arraystretch}{1.3}
\begin{array}{rcl}
\mathcal{X}\udCDdot{C}{A}{\otimes}Y\,(c) & = & \mathcal{X}(\qquest,c)\dCdot{\smallqquest\in A}{\otimes}Y(\qquest) \\
& = & \coeq\Bigg(\;\coprod_{(b\stackrel{\!\alpha}{\rightarrow}b')\in
\cat{A}}\mathcal{X}(b',c)\odot Y(b)
\;\;{\renewcommand{\arraystretch}{0.5}
\begin{array}{@{}c@{}}
\longrightarrow \\
\longrightarrow \\
\end{array}
\renewcommand{\arraystretch}{1}}\;\;
\coprod_{a\in\cat{A}}\mathcal{X}(a,c)\odot Y(a)\;\Bigg)\,. \\
\end{array}
$$
We often designate $\mathcal{X}(\qquest,?)\udCDdot{C}{\smallqquest\in A}{\otimes}
Y(\qquest)$ simply by $\mathcal{X}\udCDdot{C}{A}{\otimes}Y$; this should cause no confusion\,: the
tensor product is performed over the \emph{contravariant} variable of $\mathcal{X}$.
The coupling obtained in the other order, \ie first Lemma~\ref{odot-induced1-lem} and then~\ref{odot-induced2-lem}, is
isomorphic to the above one up to precomposition with the obvious switch functor $\sigma_{\cat{C},\cat{A}\op}\,
\times\cat{S}_{2}^{\cat{A}}$ (see~\ref{switch-not}). Indeed, unfolding the definitions gives back the same formula.

\begin{Ex}
\label{bi-tensor-ex}
Let $\cat{S}$ be a cofibrantly generated simplicial model category, and let $\boxdot$ and $\odot$
be the couplings of Examples \ref{coupling-sets-ex} and \ref{coupling-ex}. Given two small categories
$\cat{A}$ and $\cat{C}$, we get two bi-tensor products
$$
\udCDdot{C}{A}{\boxtimes}\colon\Sets^{\cat{A}\op\times\cat{C}}\times\cat{S}^{\cat{A}}
\longrightarrow\cat{S}^{\cat{C}}\qquad\mbox{and}\qquad\udCDdot{C}{A}{\otimes}\colon
\sSets^{\cat{A}\op\times\cat{C}}\times\cat{S}^{\cat{A}}\longrightarrow\cat{S}^{\cat{C}}\,.
$$
Given $Y\in\cat{S}^{\cat{A}}$ and $Z\in\cat{S}^{\cat{C}}$, we get for $\mathcal{X}\in\Sets^
{\cat{A}\op\times\cat{C}}$, a natural bijection
$$
\mor_{\cat{S}^{\cat{C}}}\!\Big(\mathcal{X}\udCDdot{C}{A}{\boxtimes}Y,Z\Big)\cong\mor_{\Sets^{\cat
{A}\op\times\cat{C}}}\!\Big(\mathcal{X},\rmap{\smalludCDdot{C}{A}{\boxtimes}}(Y,Z)\Big)
$$
and for $\mathcal{X}'\in\sSets^{\cat{A}\op\times\cat{C}}$, a natural
bijection
$$
\mor_{\cat{S}^{\cat{C}}}\!\Big(\mathcal{X}'\udCDdot{C}{A}{\otimes}Y,Z\Big)\cong\mor_{\sSets^{\cat
{A}\op\times\cat{C}}}\!\Big(\mathcal{X}',\rmap{\smalludCDdot{C}{A}{\otimes}}(Y,Z)\Big)\,.
$$
Furthermore, the right mapping functors are explicitly given by
$$
\rmap{\smalludCDdot{C}{A}{\boxtimes}}(Y,Z)(\qquest,?)=\mor_{\cat{S}}\!\big(Y(\qquest),Z(?)\big)
$$
and
$$
\rmap{\smalludCDdot{C}{A}{\otimes}}(Y,Z)(\qquest,?)=\Map\!\big(Y(\qquest),Z(?)\big)\,.
$$
\end{Ex}

\begin{Rem}
\label{res-bi-tensor-rem}
Keep notations as in \ref{bi-tensor-propdef}. Let $\cat{D}$ be a subcategory of $\cat{C}$ and let
$$
\res_{\cat{A}\op\times\cat{D}}^{\cat{A}\op\times\cat{C}}\colon\;\cat{S}_{1}^{\cat{A}\op\times\cat{C}}
\longrightarrow\cat{S}_{1}^{\cat{A}\op\times\cat{D}}
\qquad\mbox{and}
\qquad\res_{\cat{D}}^{\cat{C}}
\colon\;\cat{S}_{3}^{\cat{C}}\longrightarrow\cat{S}_{3}^{\cat{D}}
$$
be the obvious restriction functors. Then, for $\mathcal{X}\in\cat{S}_{1}^{\cat{A}\op\times\cat{C}}$
and $Y\in\cat{S}_{2}^{\cat{A}}$, one has
$$
\res_{\cat{D}}^{\cat{C}}\big(X\udCDdot{C}{A}{\otimes}Y\big)\,=\,\big(\res_{\cat{A}\op\times\cat{D}}^{\cat{A}
\op\times\cat{C}}(X)\big) \,\udCDdot{D}{A}{\otimes}Y\quad\in\,\cat{S}_{3}^{\cat{D}}\,,
$$
as evaluation at an arbitrary object $d\in\cat{D}$ immediately shows.

\end{Rem}


\section{Couplings and model category structures}

\label{s-couplings-model}


%
We study when the functor obtained from a coupling of model categories by fixing the second variable preserves some
weak equivalences. For this purpose, the next definition turns crucial. The origin of such a concept goes back to Kan's
Homotopy Lifting Extension Theorem for categories of simplicial objects (in some category) enriched over $\sSets$,
see \cite{kan}. This was then taken as Axiom \textbf{(SM7)} for a simplicial model category by Quillen in \cite{quil}
(see Axiom (4) in Definition~\ref{simpl-mod-cat-def}).

\begin{Def}
\label{corner-map-def}
Let $\cat{S}_{1}$, $\cat{S}_{2}$ and $\cat{S}_{3}$ be model categories. Let $\odot\colon\cat{S}_{1}\times\cat{S}_{2}
\too\cat{S}_{3}$ be a
right tensorial coupling, with right mapping functor $\rmap{\odot}$. We say that $\odot$ has the
\emph{corner-map property} if the following holds\,: for every cofibration $i\colon y\longrightarrow
y'$ in $\cat{S}_{2}$ and for every fibration $p\colon z\longrightarrow z'$ in $\cat{S}_{3}$,
the induced morphism in $\cat{S}_{1}$
$$
\varphi\colon
\rmap{\odot}(y',z)\longrightarrow
\rmap{\odot}(y,z)\times_{\rmap{\odot}(y,z')}\rmap{\odot}(y',z')
$$
is a fibration, and it is a trivial fibration if either $i$ or $p$ is moreover a weak equivalence. This morphism
$\varphi$ is the ``corner-map'' to the pull-back
 induced by the morphisms $\rmap{\odot}(i,z)$ and $\rmap{\odot}(y',p)$ as follows\,:
$$
\xymatrix{
\rmap{\odot}(y',z)\ar@{-->}[rd]^-{\varphi}
\ar@/^2em/[rrd]^(.6){\rmap{\odot}(y',p)}
\ar@/_3em/[rdd]_-{\rmap{\odot}(i,z)}
\\
& \rmap{\odot}(y,z)\times_{\rmap{\odot}(y,z')}\rmap{\odot}(y',z')\ar[r]\ar[d]
& \rmap{\odot}(y',z')\ar[d]^-{\rmap{\odot}(i,z')}
\\
& \rmap{\odot}(y,z)\ar[r]_-{\rmap{\odot}(y,p)}
& \rmap{\odot}(y,z')
}
$$
\end{Def}

\begin{Ex}
\label{simpl-corner-ex}
For a simplicial model category $\cat{S}$, the right tensorial coupling $\odot$ given by the
``action'' of the category $\sSets$ on $\cat{S}$, see Example \ref{coupling-ex}, has the corner-map
property, by Axiom (4) of Definition \ref{simpl-mod-cat-def}.
\end{Ex}

\begin{Rem}
\label{flat-rem}
Let $\odot\colon\cat{S}_{1}\times\cat{S}_{2}\too\cat{S}_{3}$ be a right-tensorial coupling having the corner-map
property. Let $y\in\cat{S}_{2}$ be a \emph{cofibrant} object. Applying the above condition to the cofibration
$i\colon\varnothing\too y$, it is easy to check that the functor $\rmap{\odot}(y,-)\colon\cat{S}_{3}\too\cat{S}_{1}$
preserves (trivial) fibrations. In other words, we have a Quillen adjunction
$$
-\odot y\colon\cat{S}_{1}\adjtoo\cat{S}_{3}\noloc\rmap{\odot}(y,-)\,.
$$
In particular, the functor $-\odot y$ preserves (trivial) cofibrations. So, $-\odot y$ preserves cofibrant objects
and weak-equivalences between them, by Ken Brown's Lemma (see \cite[Lem.~1.1.12]{hov}).
\end{Rem}

\medbreak
\centerline{*\ *\ *}
\bigbreak

Of course, we want to extend the above corner-map property to categories of diagrams. The obvious statement holds
(Corollary~\ref{corner-cor} below) but is not sufficient for applications. We shall need the following improved version.
Note that the morphism $i$ of the statement is only required to be an \emph{objectwise} cofibration.

\begin{Lem}
\label{corner-lem}
Let $\cat{S}_{1}$, $\cat{S}_{2}$ and $\cat{S}_{3}$ be cofibrantly generated model categories. Let $\odot\colon
\cat{S}_{1}\times\cat{S}_{2}\too\cat{S}_{3}$ be a right tensorial coupling with the corner-map property. Let
$\cat{A}$ and $\cat{C}$ be small categories and let $\cat{D}\subset\cat{C}$ be a subcategory. Consider
$$
\rmap{\smalludCDdot{C}{A}{\otimes}}\colon
\big(\cat{S}_{2}^{\cat{A}}\big)\nnspace\op\times\cat{S}_{3}^{\cat{C}}
\longrightarrow
\cat{S}_{1}^{\cat{A}\op\times\cat{C}}
$$
the right mapping functor of the bi-tensor product
$
\udCDdot{C}{A}{\otimes}\colon
\cat{S}_{1}^{\cat{A}\op\times\cat{C}}\times\cat{S}_{2}^{\cat{A}}
\longrightarrow
\cat{S}_{3}^{\cat{C}}
$ from~\ref{bi-tensor-propdef}.
Assume that $i\colon Y\to Y'$ is an objectwise cofibration in $\cat{S}_{2}^{\cat{A}}$ and that $p\colon Z\to Z'$
in $\cat{S}_{3}^{\cat{C}}$  is a (trivial) fibration in $\USCD{{}\cat{S}_{3}}{C}{D}$.
Then the ``corner-map'' morphism
$$
\Phi\colon\rmap{\smalludCDdot{C}{A}{\otimes}}(Y',Z)\longrightarrow\rmap{\smalludCDdot{C}{A}{\otimes}}(Y,Z)
\mathop{\times}_{\rmap{\smalludCDdot{C}{A}{\otimes}}(Y,Z')}\rmap{\smalludCDdot{C}{A}{\otimes}}(Y',Z')
$$
is a (trivial) fibration in $\USCD{{}\cat{S}_{1}}{A\op\times\cat{C}\,}{\,\cat{A}\op\times\cat{D}}$.
\end{Lem}

\Prf
Recall from Proposition~\ref{bi-tensor-propdef}, where we saw that $\udCDdot{C}{A}{\otimes}$ is right tensorial, that
$$
\rmap{\smalludCDdot{C}{A}{\otimes}}(Y,Z)\,(\qquest,?):=
\rmap{\odot}\!\big(Y(\qquest),Z(?)\big)\,.
$$
Since the pull-back and the morphism $\Phi$ are defined in the diagram-category $\cat{S}_{1}^{\cat{A}\op\times\cat{C}}$,
they can be computed $\cat{A}\op\times\cat{C}$-objectwise. Let $a$ be an object of $\cat{A}$ and let $d$ be an object of
$\cat{D}\subset\cat{C}$. The ``corner-map'' $\Phi(a,d)$ is the following map\,:
$$
\kern-4.5em
\xymatrix@C=.8em{
{\kern4em}\rmap{\odot}(Y'(a),Z(d))
\ar[rd]^-{\Phi(a,d)}
\ar@/^2em/[rrrd]^(.6){\ \ \rmap{\odot}(Y'(a),p(d))}
\ar@/_4em/[rdd]_-{\rmap{\odot}(i(a),Z(d))}
\\
& \kern-2.5em
\rmap{\odot}(Y(a),Z(d))\kern-2.5em
\mathop{\times}\limits_{\rmap{\odot}(Y(a),Z'(d))^{\vphantom{Y^{Y}}}}^{\vphantom{Y}}
\kern-2.5em\rmap{\odot}(Y'(a),Z'(d))
\ar[rr]\ar[d]
&& \rmap{\odot}(Y'(a),Z'(d))\ar[d]_-{\rmap{\odot}(i(a),Z'(d))}
\\
& \rmap{\odot}(Y(a),Z(d))\ar[rr]_-{\rmap{\odot}(Y(a),p(d))}
&& \rmap{\odot}(Y(a),Z'(d))
}
$$
By assumption $i(a)\colon Y(a)\longrightarrow Y'(a)$ is a cofibration in $\cat{S}_{2}$ and $p(d)\colon Z(d)
\longrightarrow Z'(d)$ is a (trivial) fibration in $\cat{S}_{3}$.
So, $\Phi(a,d)$ indeed is a (trivial) fibration since the original coupling $\odot\colon\cat{S}_{1}\times
\cat{S}_{2}\too\cat{S}_{3}$ has the corner map property.
\qed

\begin{Cor}
\label{corner-cor}
Let $\cat{S}_{1}$, $\cat{S}_{2}$ and $\cat{S}_{3}$ be cofibrantly generated model categories. Let $\odot
\colon\cat{S}_{1}\times\cat{S}_{2}\too\cat{S}_{3}$ be a right tensorial coupling with the corner-map property.
Let $\cat{A}$ and $\cat{C}$ be small categories and let $\cat{D}\subset\cat{C}$ be a subcategory.
Then the bi-tensor product $\udCDdot{C}{A}{\otimes}$ has the cornel map property for the following model structures\,:
$$
\udCDdot{C}{A}{\otimes}\colon\;
\USCD{{}\cat{S}_{1}}{A\op\times\cat{C}}{A\op\times\cat{D}}
\times
\USC{{}\cat{S}_{2}}{A}
\too
\USCD{{}\cat{S}_{3}}{C}{D}\,.
$$
\end{Cor}

\Prf
This is a direct consequence of Lemma~\ref{corner-lem} since a cofibration $i\colon Y\to Y'$ in $\USC{{}\cat{S}_{2}}{A}$
is in particular $\cat{A}$-objectwise a cofibration (see Remark~\ref{objwise-cof-rem}).
\qed

\goodbreak \bigbreak
\centerline{*\ *\ *}
\bigbreak

Our next goal is to single out some situations where taking the coupling with a given object preserves weak equivalence.

\begin{Thm}
\label{flat-thm}
Let $\cat{S}_{1}$, $\cat{S}_{2}$ and $\cat{S}_{3}$ be cofibrantly generated model categories and let
$\odot\colon\cat{S}_{1}\times\cat{S}_{2}\too\cat{S}_{3}$ be a right tensorial coupling having the corner-map property.
Let $\cat{A}$ and $\cat{C}$ be small categories, and $\cat{D}$
a subcategory of $\cat{C}$.
Let $Y\in\cat{S}_{2}^{\cat{A}}$ be $\cat{A}$-objectwise cofibrant. Then, there
is a Quillen adjunction
$$
(-)\udCDdot{C}{A}{\otimes}Y\colon\quad
\UUU{\cat{S}_{1}}{\cat{A}\op\times\cat{C}}{\cat{A}\op\times\cat{D}}
\,\adjtoo\,
\USCD{S_{3}}{C}{D}
\quad\noloc\rmap{\smalludCDdot{C}{A}{\otimes}}(Y,-)\,.
$$
In particular, in this situation, the functor $(-)\udCDdot{C}{A}{\otimes}Y$ takes $\cat{A}\op\times\cat{D}$-weak
equivalences between cofibrant objects in $\UUU{\cat{S}_{1}}{\cat{A}\op\times\cat{C}}{\cat{A}\op\times\cat{D}}$
to $\cat{D}$-weak equivalences between cofibrant objects in $\USCD{S_{3}}{C}{D}$.
\end{Thm}

\Prf
It suffices to show that
$\rmap{\smalludCDdot{C}{A}{\otimes}}(Y,-)\colon\USCD{S_{3}}{C}{D}\too\UUU{\cat{S}_{1}}{\cat{A}\op\times\cat{C}}{\cat{A}
\op\times\cat{D}}$ preserves fibrations and trivial fibrations, which are tested objectwise in both the source and target
model categories. Let $p\colon Z\longrightarrow Z'$ be a (trivial) $\cat{D}$-fibration in $\cat{S}_{3}^{\cat{C}}$.
We conclude by Lemma~\ref{corner-lem} applied to $i\colon\varnothing\too Y$ and to our $p\colon Z\too Z'$.
\qed

\medbreak

We leave it to the reader to unfold the obvious corollaries of the results of the present section when $\cat{A}$ is reduced to a
point or when $\cat{D}=\cat{C}$. We only mention for further quotation the following consequence of Theorem~\ref{flat-thm}
with $\cat{D}=\cat{C}=\{*\}$.

\begin{Cor}
\label{flat-cor}
Let $\cat{S}_{1}$, $\cat{S}_{2}$ and $\cat{S}_{3}$ be cofibrantly generated model categories and let
$\odot\colon\cat{S}_{1}\times\cat{S}_{2}\too\cat{S}_{3}$ be a right tensorial coupling having the corner-map property.
Let $\cat{A}$ be a small category.
Let $Y\in\cat{S}_{2}^{\cat{A}}$ be $\cat{A}$-objectwise cofibrant. Let $\xi\colon X\too X'$ in $\cat{S}_{1}^{\cat{A}\op}$
be a weak-equivalence between cofibrant objects in $\USC{S_{1}}{A\op}$. Then the morphism $\xi\dCdot{A}{\otimes}Y\colon
X\dCdot{A}{\otimes}Y\too X'\dCdot{A}{\otimes}Y$ is a weak equivalence in $\cat{S}_{3}$.
\qed
\end{Cor}


\section{Cofibrant approximations in $\USCD{S}{C}{D}$}

\label{s-cof-approx}


The goal in this section is to establish a general technique to construct functorial cofibrant approximation in the model
category $\USCD{S}{C}{D}$. Explicit examples will be presented in Section \ref{s-var-cof-approx}. The reader who prefers
to focus on the model category $\USC{S}{C}$ may well take $\cat{D}=\cat{C}$ everywhere in the present section, although
in this case, some proofs look more transparent when distinguishing two different occurrences of $\cat{C}$, calling one
of them $\cat{D}$ and the other one $\cat{C}$ as we do here.

\begin{Not}
\label{mor-eta-not}
Let $\cat{C}$ be a small category. We fix the following notations for the rest of the paper.
\begin{itemize}
\item [(i)] Let $\cat{S}$ be a model category with cofibrant replacement denoted by $(Q_{\cat{S}},\xi^{\cat{S}})$.
For a diagram $X\in\cat{S}^{\cat{C}}$, we denote by $qX$ the diagram
$Q_{\cat{S}}\circ X\in\cat{S}^{\cat{C}}$ obtained by applying the cofibrant replacement
$\cat{C}$-objectwise to $X$, that is,
$$
\itemspace\xymatrix{
qX\colon\cat{C} \ar[r]^-{X} & \cat{S} \ar[r]^-{Q_{\cat{S}}} & \cat{S}
},
\qquad qX(c)=Q_{\cat{S}}\big(X(c)\big)\,.
$$
We let $\eta_{X}\colon qX\longrightarrow X$ be the morphism given, for $c\in\cat{C}$, by
$\eta_{X}(c):=\xi^{\cat{S}}_{X(c)}$.
\item [(ii)]  For a subcategory $\cat{D}$ of $\cat{C}$, we define the diagram of sets
$$
\FDC{D}{C}\colon\;
\cat{D}\op\times\cat{C}\longrightarrow\Sets,
\qquad
(d,c) \longmapsto\mor_{\cat{C}}(d,c)\,.
$$
Composing with the usual constantification functor $\Cst{\Sets}{\DDelta\!\op}\colon\Sets
\longrightarrow\sSets$, we obtain the functor $\FFDC{D}{C}:=\underline{\FDC{D}{C}}$, that is,
$$
\itemspace
\FFDC{D}{C}\colon\;
\cat{D}\op\times\cat{C}\longrightarrow\sSets,
\qquad(d,c)\longmapsto\underline{\mor_{\cat{C}}(d,c)}\,.
$$
When $\cat{D}=\cat{C}$ we shall abbreviate $\FFDC{C}{C}$ by $\FFC{C}$ but we still write $\FDC{C}{C}$ to stress
the difference between this diagram of sets and  the mere set $\mor_{\cat{C}}$.
\end{itemize}
\end{Not}

Here is the first result regarding the diagram $\FFDC{D}{C}$, giving it some flavour.

\begin{Lem}
\label{morCD-ind-lem}
Let $\cat{D}$ be a full subcategory of a small category $\cat{C}$. Then, there are canonical isomorphisms
$$
\ind_{\cat{D}\op\times\cat{D}}^{\cat{D}\op\times\cat{C}}\!\big(\FDC{D}{D}\big)\cong\FDC{D}{C}
\ \in\Sets^{\cat{D}\op\times\cat{C}}
$$
and
$$
\ind_{\cat{D}\op\times\cat{D}}^{\cat{D}\op\times\cat{C}}\big(\FFC{D}\big)
\cong\FFDC{D}{C}\
\in\sSets^{\cat{D}\op\times\cat{C}}\,.
$$
From now on, for simplicity, we consider these identifications as equalities.
\end{Lem}

\Prf
For an object $d\in\cat{D}\op$, it suffices to produce in $\Sets^{\cat{C}}$ an isomorphism
$
\ind_{\cat{D}}^{\cat{C}}\!\big(\mor_{\cat{D}}(d,?)\big)\too
\FDC{D}{C}(d,?)=\mor_{\cat{C}}(d,?)
$
that is natural in~$d$. The induction functor being a left Kan extension
(see \cite{MacLane} or~\cite[Def.\,\ref{LKan-def}]{bamaI}), we have for all $c\in\cat{C}$
$$
\Big(\ind_{\cat{D}}^{\cat{C}}\!\big(\mor_{\cat{D}}(d,?)\big)\Big)(c)\cong\colim_{\cat{D}\smallcomma c}\,
\mor_{\cat{D}}(d,d')=\colim_{\cat{D}\smallcomma c}\,\mor_{\cat{C}}(d,d')\cong\mor_{\cat{C}}(d,c)\,,
$$
where the colimits are taken over $(d',d'\stackrel{\alpha}{\rightarrow}c)\in\cat{D}\comma c$. The equality in
the middle follows from fullness. The last isomorphism is easy. The second part follows.
\qed

\begin{Lem}
\label{lambda-iso-lem}
Let $\cat{S}$ be a simplicial model category and $\cat{D}$ a full subcategory of a small category $\cat{C}$.
Let $\udCDdot{C}{D}{\boxtimes}$ and $\udCDdot{C}{D}{\otimes}$ be the bi-tensor products induced
by the couplings $\boxdot$ of \ref{coupling-sets-ex} and $\odot$ of \ref{coupling-ex}, respectively.
Then, for every diagram $X\in\cat{S}^{\cat{D}}$, there are canonical and natural isomorphisms
$$
\xymatrix{
\lambda_{X}^{\cat{D},\cat{C}}\colon
\quad\FFDC{D}{C}(\qquest,?)\,\udCDdot{C}{D}{\otimes}\,X(\qquest)
\ar[r]^-{\cong}
& \FDC{D}{C}(\qquest,?)\,\udCDdot{C}{D}{\boxtimes}\,X(\qquest)
\ar[r]^-{\cong}
& \ind_{\cat{D}}^{\cat{C}}X\,(?)\,.
}
$$
\end{Lem}

\Prf
The first isomorphism follows from Lemma \ref{compatibility-lem}. So, it suffices to construct the second isomorphism, for
which it is enough to check that the functor
$
\FDC{D}{C}\udCDdot{C}{D}{\boxtimes}-\,\colon\,\cat{S}^{\cat{D}}\too\cat{S}^{\cat{C}}
$
is a left adjoint of $\res_{\cat{D}}^{\cat{C}}$. Consider two diagrams $Y\in\cat
{S}^{\cat{D}}$ and $Z\in\cat{S}^{\cat{C}}$.  Using the adjunctions and the explicit formulas of~\ref{bi-tensor-ex} as well as
Lemma~\ref{morCD-ind-lem} (at the second step), we have natural bijections
$$
\arraycolsep1pt
\renewcommand{\arraystretch}{1.4}
\begin{array}{rcl}
\mor_{\cat{S}^{\cat{C}}}\!\Big(\FDC{D}{C}\udCDdot{C}{D}{\boxtimes}\,Y\,,\,Z\Big) & \cong &
\mor_{\Sets^{\cat{D}\op\times\cat{C}}}\!\Big(\FDC{D}{C}\,,\rmap{\smalludCDdot{C}{D}{\boxtimes}}(Y,Z)\Big)
\\
& = & \mor_{\Sets^{\cat{D}\op\times\cat{C}}}\!\Big(\ind_{\cat{D}\op\times\cat{D}}^{\cat{D}\op\times\cat{C}}\!
\big(\FDC{D}{D}\big),\rmap{\smalludCDdot{C}{D}{\boxtimes}}(Y,Z)\Big) \\
& \cong & \mor_{\Sets^{\cat{D}\op\times\cat{D}}}\!\Big(\FDC{D}{D}\,,\res_{\cat{D}\op\times\cat{D}}^{\cat{D}\op\times\cat{C}}
\!\big(\rmap{\smalludCDdot{C}{D}{\boxtimes}}(Y,Z)\big)\Big) \\
& = & \mor_{\Sets^{\cat{D}\op\times\cat{D}}}\!
\Big(\mor_{\cat{D}}(\qquest,?)\ ,\ \mor_{\cat{S}}\!\big(Y(\qquest),\res_{\cat{D}}^{\cat{C}}Z(?)\big)\Big) \\
& \cong & \mor_{\cat{S}^{\cat{D}}}\!\big(Y,\res_{\cat{D}}^{\cat{C}}Z\big)\,.
\\
\end{array}
$$
The final bijection is an easy exercise of Yoneda style.
\qed

\goodbreak
\medbreak

For the next statements, we adopt the notations of~\ref{mor-eta-not} and~\ref{lambda-iso-lem}.

\begin{Thm}
\label{cof-approx-USCD-thm}
Let $\cat{D}$ be a full subcategory of a small category $\cat{C}$, and $\cat{S}$ a cofibrantly generated
simplicial model category with the ``action'' of $\sSets$ denoted by $\odot$ (see \ref{coupling-ex}).
Denote by $\udCDdot{C}{D}{\otimes}\colon
\sSets^{\cat{D}\op\times\cat{C}}
\times\cat{S}^{\cat{D}}
\longrightarrow\cat{S}^{\cat{C}}$
the bi-tensor product induced by $\odot$
(see \ref{bi-tensor-propdef} and \ref{bi-tensor-ex}).
In $\UUU{\sSets}{\cat{D}\op\times\cat{C}\,}{\,\cat{D}\op\times\cat{D}}$, fix a cofibrant approximation
$$
\vartheta\colon\bbE\longrightarrow\FFDC{D}{C}
$$
of the diagram $\FFDC{D}{C}$ given by $\FFDC{D}{C}(d,c)=\underline{\mor_{\cat{C}}(d,c)}$ for all $d\in\cat{D}$ and $c\in\cat{C}$.
More precisely, $\bbE$ is an arbitrary cofibrant object in this model category and $\vartheta$ is a weak equivalence $\cat{D}\op
\times\cat{D}$-objectwise. For every
diagram $X\in\cat{S}^{\cat{C}}$, we define
$$
\mQCD{C}{D}X:=\bbE\,\udCDdot{C}{D}{\otimes}\res_{\cat{D}}^{\cat{C}}qX
$$
and we let $\xiCDX{C}{D}{X}\colon\mQCD{C}{D}X\too X$ be given by the composition
$$
\xymatrix @C=3.2em @R=1.6em{
{\bbE}\,\udCDdot{C}{D}{\otimes}\res_{\cat{D}}^{\cat{C}}qX
\ar[r]^-{\vartheta\smalludCDdot{C}{D}{\otimes}\res\nnspace q\nnspace X}
\ar@/_2.6em/[rrrr]|{\ \xiCDX{C}{D}{X}\ }
& \FFDC{D}{C}\ \udCDdot{C}{D}{\otimes}\,
\res_{\cat{D}}^{\cat{C}}qX
\ar[r]^-{\lambda_{\res\nnspace q\nnspace X}^{\cat{D},\cat{C}}}_-{\cong}
&\ind_{\cat{D}}^{\cat{C}}\res_{\cat{D}}^{\cat{C}}qX
\ar[r]^-{\epsilon_{q\nnspace X}}
& qX
\ar[r]^(.45){\eta_{X}}_(.45){\sim}
& X
}
$$
where $\epsilon_{q\nnspace X}$ is the counit, at $qX$, of the adjunction $\big(\ind_{\cat{D}}^{\cat{C}}\,,\,
\res_{\cat{D}}^{\cat{C}}\big)$. Then, the pair $(\mQCD{C}{D},\xiCD{C}{D})$ is a functorial cofibrant
approximation in the model category $\USCD{S}{C}{D}$.
Moreover, the morphism $\vartheta\smalludCDdot{C}{D}{\otimes}qX$ is a $\cat{D}$-weak equivalence,
$\epsilon_{qX}$ is a $\cat{D}$-isomorphism, and $\eta_{X}$ is a $\cat{C}$-weak equivalence.
\end{Thm}

After a couple of remarks and immediate corollaries, the proof will be given at the end of this section, starting with
Lemma~\ref{EC-objwise-cof-lem}.

\begin{Rem}
\label{cof-approx-USCD-rem}
There is no need for $\bbE$ and $\vartheta$ to be part of a functorial cofibrant approximation\,: one only needs
to know how to cofibrantly approximate
the \emph{single} diagram $\FFDC{D}{C}$ in the model category $\UUU{\sSets}
{\cat{D}\op\times\cat{C}}{\cat{D}\op\times\cat{D}}$. In fact, it follows from Remark~\ref{res-ind-Q-adj-rem}
that there is a Quillen adjunction
$$
\ind_{\cat{D}\op\times\cat{D}}^{\cat{D}\op\times\cat{C}}\colon\quad\UUSC{\sSets}{\cat{D}\op\times\cat{D}}\adjtoo
\UUU{\sSets}{\cat{D}\op\times\cat{C}}{\cat{D}\op\times\cat{D}}\quad\noloc\res_{\cat{D}\op\times\cat{D}}^{\cat{D}
\op\times\cat{C}}\,.
$$
In particular, it is enough to select a cofibrant approximation $\varrho\colon\bbE_{0}\longrightarrow\FFC{D}$
in the model category $\UUSC{\sSets}{\cat{D}\op\times\cat{D}}$ and to induce it up to get
$$
\underbrace{\ind_{\cat{D}\op\times\cat{D}}^{\cat{D}\op\times\cat{C}}(\varrho)}_{=:\vartheta}\colon
\underbrace{\ind_{\cat{D}\op\times\cat{D}}^{\cat{D}\op\times\cat{C}}(\bbE_{0})}_{=:\bbE}\longrightarrow
\ind_{\cat{D}\op\times\cat{D}}^{\cat{D}\op\times\cat{C}}\!\big(\FFC{D}\big)
=\FFDC{D}{C}\,,
$$
where the latter identification is given by Lemma \ref{morCD-ind-lem}; then, $\vartheta$ and $\bbE$
fulfill the required properties needed in Theorem \ref{cof-approx-USCD-thm}. The point here is that
the model category $\UUSC{\sSets}{\cat{D}\op\times\cat{D}}$ is {\sl à la} Dwyer-Kan-Heller-Dugger as
considered in \cite{dwyerkan1a,hell,dugg}, that is, with $\sSets$ as category of ``values'' and without
need of a relative model category structure in the sense of \ref{USCD-propdef}. See also Remark~\ref{various-rem} below.
\end{Rem}

\begin{Rem}
\label{choice-QS-rem}
The only property of the morphism $\eta_{X}\colon qX\longrightarrow X$, where $X\in\cat{S}^{\cat{C}}$,
that is needed here is that it is a $\cat{D}$-weak equivalence of $\cat{C}$-indexed diagrams, depending functorially
on $X$, with a $\cat{D}$-objectwise cofibrant source $qX$. Therefore, in this section, one could have replaced
everywhere the cofibrant replacement $(Q_{\cat{S}},\xi^{\cat{S}})$ in the model category $\cat{S}$ by any
functorial cofibrant approximation $(\mathcal{Q}_{\cat{S}},\zeta^{\cat{S}})$; the functoriality of $\zeta^{\cat{S}}$
guarantees that $qX$ really is a functor. One could even merely chose an arbitrary $\cat{D}$-objectwise
cofibrant diagram $qX$ in $\cat{S}^{\cat{C}}$ and an arbitrary $\cat{D}$-weak equivalence $\eta_{X}\colon qX
\longrightarrow X$ in $\cat{S}^{\cat{C}}$ at the only cost of dropping the functoriality statement in
Theorem~\ref{cof-approx-USCD-thm}. In particular, if $X$ is already known to be
$\cat{D}$-objectwise cofibrant, it suffices to take $qX=X$ and $\eta_{X}=\id_{X}$. This gives the immediate\,:
\end{Rem}

\begin{Cor}
\label{X-objcof-cor}
Let $\cat{D}$ be a full subcategory of a small category $\cat{C}$, and $\cat{S}$ a cofibrantly generated simplicial model
category. Fix as above a cofibrant approximation
$
\vartheta\colon\bbE\longrightarrow\FFDC{D}{C}
$
in the model category
$\UUU{\sSets}{\cat{D}\op\times\cat{C}\,}{\,\cat{D}\op\times\cat{D}}$. For every
\emph{$\cat{C}$-objectwise cofibrant} diagram $X\in\cat{S}^{\cat{C}}$, we define
$$
\nQX{X}:=\bbE\,\udCDdot{C}{D}{\otimes}\res_{\cat{D}}^{\cat{C}}X
$$
and we let $\nzetaX{X}\colon\nQX{X}\too X$ be given by the composition
$$
\xymatrix @C=3.2em @R=1.6em{
{\bbE}\,\udCDdot{C}{D}{\otimes}\res_{\cat{D}}^{\cat{C}}X
\ar[r]^-{\vartheta\smalludCDdot{C}{D}{\otimes}\res\nnspace X}
\ar@/_2em/[rrr]|{\ \nzetaX{X}\ }
& \FFDC{D}{C}\ \udCDdot{C}{D}{\otimes}\,
\res_{\cat{D}}^{\cat{C}}X
\ar[r]^-{\lambda_{\res\nnspace X}^{\cat{D},\cat{C}}}_-{\cong}
&\ind_{\cat{D}}^{\cat{C}}\res_{\cat{D}}^{\cat{C}}X
\ar[r]^-{\epsilon_{X}}
& X\,.
}
$$
Then, $\nzetaX{X}\colon\nQX{X}\too X$ is a cofibrant
approximation in $\USCD{S}{C}{D}$.
\qed
\end{Cor}

We state the absolute case $\cat{D}=\cat{C}$ of Theorem \ref{cof-approx-USCD-thm}, using that
$\res_{\cat{C}}^{\cat{C}}=\id_{\cat{C}}=\ind_{\cat{C}}^{\cat{C}}$.

\begin{Cor}
\label{cof-approx-USC-cor}
Let $\cat{C}$ be a small category and $\cat{S}$ a cofibrantly generated simplicial model category. Fix a cofibrant approximation
$$
\vartheta\colon\bbE\longrightarrow\FFC{C}
$$
of the diagram $\FFC{C}(\qquest,?)=\underline{\mor_{\cat{C}}(\qquest,?)}$ in the model category $\UUSC{\sSets}{\cat{C}
\op\times\cat{C}}$; more precisely, $\bbE$ is a cofibrant object in the latter model category and
$\vartheta$ is a weak equivalence $\cat{C}\op\times\cat{C}$-objectwise. For every diagram $X\in\cat{S}^{\cat{C}}$, we define
$$
\mQC{C}X:=\bbE\,\udCDdot{C}{C}{\otimes}qX
$$
and we let $\xiCX{C}{X}$ be given by the composition
$$
\xymatrix @C=4.2em{
{\bbE}\udCDdot{C}{C}{\otimes}qX \ar[r]^-{\vartheta\smalludCDdot{C}{C}{\otimes}q\nnspace X}
\ar@/_2em/[rrr]|{\ \xiCX{C}{X}\ } & {\FFC{C}}\,\udCDdot{C}{C}
{\otimes}\,qX\ar[r]^-{\lambda_{q\nnspace X}^{\cat{C},\cat{C}}}_-{\cong} & qX \ar[r]^-{\eta_
{X}}_-{\sim} & X\,.
}
$$
Then, the pair $(\mQC{C},\xiC{C})$ is a functorial cofibrant approximation in $\USC{S}{C}$.
\qed
\end{Cor}

\medbreak
\centerline{*\ *\ *}
\medbreak

We establish now a few preparatory results for the proof of Theorem \ref{cof-approx-USCD-thm}.

\begin{Lem}
\label{EC-objwise-cof-lem}
Let $\cat{D}$ be a subcategory of a small category $\cat{C}$, and let $\bbE$ be as in \ref{cof-approx-USCD-thm}.
Consider two objects $d\in\cat{D}$ and $c\in\cat{C}$. Then, the diagrams $\bbE(\qquest,c)$ and $\bbE(d,?)$ are
cofibrant in $\UUSC{\sSets}{\cat{D}\op}$ and $\UUSC{\sSets}{\cat{C}}$ respectively.
\end{Lem}

\Prf
This is merely a special case of Remark \ref{CD-obj-wise-cof-rem}.
\qed

\begin{Lem}
\label{mor-objwise-cof-lem}
Let $\cat{C}$ be a small category, $\cat{D}$ be a subcategory of $\cat{C}$ and $d$ an object of $\cat{D}$. Then, the diagram
$\FFC{C}(\qquest,d)\in\sSets^{\cat{C}\op}$ is cofibrant in
$\USCD{{}\sSets}{C\op}{D\op}$.
\end{Lem}

\Prf
First, recall that $\FFC{C}(\qquest,d)=\underline{\mor_{\cat{C}}
(\qquest,d)}$. We will show that $\varnothing\too\FFC{C}(\qquest,d)$ has
the left lifting property with respect to trivial fibrations in
$\USCD{{}\sSets}{C\op}{D\op}$.
Take $p\colon A\longrightarrow B$ a trivial fibration, \ie a
trivial fibration $\cat{D}$-objectwise. We only have to show that the map of sets
$$
\mor_{\cat{T}}\!\big({\FFC{C}(\qquest,d)},p(\qquest)\big)\colon
\mor_{\cat{T}}\!\big({\FFC{C}(\qquest,d)},A(\qquest)\big)\longrightarrow
\mor_{\cat{T}}\!\big({\FFC{C}(\qquest,d)},B(\qquest)\big)
$$
is surjective, where we abbreviate $\cat{T}:=\sSets^{\cat{C}\op}$. We have natural bijections
$$
\mor_{\sSets^{\cat{C}\op}}\!\Big(\FFC{C}(\qquest,d),A(\qquest)\Big)
\cong
\mor_{\Sets^{\cat{C}\op}}\!\big(\mor_{\cat{C}}(\qquest,d),A(\qquest)_{0}\big)
\cong A(d)_{0}\,,
$$
using the adjunction~\ref{Cst-level-zero-rem} and the Yoneda Lemma. Doing the same with $B(\qquest)$
and with $p(\qquest)$, we are left to verify that the map
$
p(d)_{0}\colon A(d)_{0}\longrightarrow B(d)_{0}
$
is surjective in $\Sets$. This is immediate from the fact that $p(d)$ is a trivial Kan fibration in $\sSets$.
(This is well-known and follows for instance from the right lifting property with respect to the cofibration $\varnothing
\hookrightarrow\Delta^{0}$\,.)
\qed

\begin{Lem}
\label{weq-OK-lem}
Let $\cat{D}$ be a subcategory of a small category $\cat{C}$, and $\cat{S}$ a cofibrantly generated
simplicial model category. Let $\bbE$ and $\vartheta$ be as in \ref{cof-approx-USCD-thm}. Suppose
given $Y\in\cat{S}^{\cat{C}}$, a diagram that is $\cat{D}$-objectwise cofibrant. Then, the morphism
$$
\vartheta\udCDdot{C}{D}{\otimes}Y\colon\;\bbE\,\udCDdot{C}{D}{\otimes}Y\longrightarrow
\FFDC{D}{C}\,\udCDdot{C}{D}{\otimes}\,Y
$$
is a cofibrant approximation in the model category $\USCD{S}{C}{D}$, that is, the diagram
$\bbE\,\udCDdot{C}{D}{\otimes}Y$ is $\cat{D}$-cofibrant and the morphism $\vartheta\udCDdot{C}{D}
{\otimes}Y$ is a $\cat{D}$-weak equivalence.
\end{Lem}

\Prf
Since $\bbE$ is $\cat{D}\op\times\cat{D}$-cofibrant and $Y$ is $\cat{D}$-objectwise cofibrant, the
diagram $\bbE\,\udCDdot{C}{D} {\otimes}Y$ is cofibrant by Theorem~\ref{flat-thm} applied with
$\cat{A}:=\cat{D}$. It only remains to show that $\vartheta\udCDdot{C}{D}{\otimes}Y$ is a $\cat{D}$-weak
equivalence. For all $d\in\cat{D}$, the morphism $\big(\vartheta\udCDdot{C}{D}{\otimes}Y\big)\,(d)$ is equal to
$\vartheta(\qquest,d)\dCdot{D}{\otimes}Y(\qquest)\colon \bbE\,(\qquest,d)\dCdot{D}{\otimes}Y(\qquest)\too{\FFDC{D}
{C}(\qquest,d)}\dCdot{D}{\otimes}Y(\qquest)$.
This is a weak-equivalence in $\cat{S}$ by Corollary~\ref{flat-cor} for $\cat{A}:=\cat{D}$ and $\xi(\qquest):=
\vartheta(\qquest,d)$. To apply \ref{flat-cor}, we need the following facts\,: $Y$ is objectwise cofibrant by
hypothesis; the diagram $\bbE\,(\qquest,d)$ is cofibrant in $\USC{{}\sSets}{D\op}$ by Lemma~\ref{EC-objwise-cof-lem};
one has ${\FFDC{D}{C}(\qquest,d)}=\FFC{D}(\qquest,d)$ in $\sSets^{\cat{D}\op}$ and hence ${\FFDC{D}{C}(\qquest,d)}$ is
cofibrant in $\USC{{}\sSets}{D\op}$ by the absolute version of Lemma~\ref{mor-objwise-cof-lem}; and finally $\vartheta
(\qquest,d)\colon\bbE\,(\qquest,d)\too{\FFDC{D}{C}(\qquest,d)}$
is a weak-equivalence in $\USC{{}\sSets}{D\op}$ by definition of~$\vartheta$.
\qed

\medskip

\PrfOf{Theorem \ref{cof-approx-USCD-thm}}
By Lemma~\ref{weq-OK-lem} applied to $Y:=\res_{\cat{D}}^{\cat{C}}qX$, which is $\cat{D}$-objectwise cofibrant,
the diagram $\bbE\,\udCDdot{C}{D}{\otimes}\res_{\cat{D}}^{\cat{C}}qX$
is $\cat{D}$-cofibrant and $\vartheta\udCDdot{C}{D}{\otimes}qX$ is a $\cat{D}$-weak equivalence. Since the
functor $\res_{\cat{D}}^{\cat{C}}\circ\ind_{\cat{D}}^{\cat{C}}$ is naturally equivalent to $\id_{\cat{D}}$, the
morphism $\res_{\cat{D}}^{\cat{C}}(\epsilon_{q\nnspace X})$ is an isomorphism, that is, $\epsilon_
{q\nnspace X}$ is a $\cat{D}$-isomorphism. Since $\eta_{X}$ is a $\cat{D}$-weak equivalence (and even a
$\cat{C}$-weak equivalence), $\xiCDX{C}{D}{X}$ is a $\cat{D}$-weak equivalence and the result follows.
\qed


\section{Explicit examples of cofibrant approximations in $\USCD{S}{C}{D}$}

\label{s-var-cof-approx}


In the present section, we provide an explicit cofibrant approximations of
the diagram $\FFDC{D}{C}$ in $\USCD{{}\sSets}{D\op\times\cat{C}}{D\op\times\cat{D}}$. With Theorem~\ref{cof-approx-USCD-thm},
this produces an explicit functorial cofibrant approximations in $\USCD{S}{C}{D}$.

\medskip

For this section, we fix a cofibrantly generated simplicial model category
$\cat{S}$.

\medskip

We need some notations. To start with, given a small category
$\cat{C}$, we denote by $B\cat{C}=B_{\bullet}\cat{C}\in\sSets$ its \emph{nerve}, whose
realization $|B\cat{C}|\in\Top$ is the usual \emph{classifying space} of $\cat{C}$. Here,
we follow Segal's modern definition of the nerve in \cite{segal}, see also Quillen~\cite{quil2} (and
not Bousfield-Kan's old definition in \cite{boukan}, where their $B\cat{C}$ is our $B(\cat{C}\op)$;
note however that $|B\cat{C}|$ and $|B(\cat{C}\op)|$ are canonically homeomorphic).

\begin{Not}
\label{commaC-not}
Let $\cat{D}$ be a subcategory of a small category $\cat{C}$.
\begin{itemize}
\item [(i)] Given two objects $d\in\cat{D}$ and $c\in\cat{C}$, we let $d\comma\cat{D}\commaC{C}c$ be the
double-comma category with
$$
\itemspace\big\{d\stackrel{\alpha}{\longrightarrow}d_{0}\stackrel{\gamma}{\longrightarrow}c\,\big|\,
d_{0}\in\obj(\cat{D}),\,\alpha\in\arr(\cat{D})\;\,\mbox{and}\;\,\gamma\in\arr(\cat{C})\big\}
$$
as set of objects, with the commutative diagrams of the form
$$
\itemspace\xymatrix @R=1.5em{
d \ar[r]^-{\alpha} \ar@{=}[d] & d_{0} \ar[r]^-{\gamma} \ar[d]_{\beta} & d \ar@{=}[d] \\
d \ar[r]^-{\alpha'} & d_{0}' \ar[r]^-{\gamma'} & d \\
}
$$
with $\beta\in\arr(\cat{D})$, as morphisms, and with the obvious concatenation of diagrams as composition.
The comma category $\cat{D}\commaC{C}c$ is defined similarly.
\item [(ii)] In the diagram-category $\sSets^{\cat{D}\op\times\cat{C}}$, consider the diagram $\EEDC{D}{C}=\EEDC{D}{C}
(\qquest,?)$ defined by
$$
\itemspace \EEDC{D}{C}:=B(\qquest\comma\cat{D}\commaC{C}?)\op\colon\;
\cat{D}\op\times\cat{C}
\longrightarrow\sSets,\qquad(d,c)\longmapsto B(d\comma\cat{D}\commaC{C}c)\op
$$
and the morphism $\vthetaDC{D}{C}\colon\EEDC{D}{C}\longrightarrow
\FFDC{D}{C}$ given at level $n\in\bbN$ by the obvious composition, namely,
for $d\in\cat{D}$ and $c\in\cat{C}$,
$$
\itemspace
\vcenter{\xymatrix @R=1.5em{
d \ar[r]^-{\smash[t]{\alpha_{0}}} \ar@{=}[d] & d_{0} \ar[r]^{\smash[t]{\beta_{0}}} \ar@{<-}[d]^-{\!\gamma_{1}} &
c \ar@{=}[d] \\
d \ar[r]^-{\alpha_{1}} \ar@{=}[d] & d_{1} \ar[r]^{\beta_{1}}
\ar@{<-}[d]^-{\!\gamma_{2}} & c \ar@{=}[d] \\
{\raisebox{.6em}{\vdots}} \ar@{=}[d] & {\raisebox{.6em}{\vdots}} \ar@{<-}[d]^-{\!\gamma_{n}} &
{\raisebox{.6em}{\vdots}} \ar@{=}[d] \\
d \ar[r]^-{\alpha_{n}} & d_{n} \ar[r]^{\beta_{n}} & c \\
}}
\qquad
\vcenter{\xymatrix @C=4.5em{
{\phantom{d}} \ar@{|->}[r]^-{\smash[t]{\vthetaDC{D}{C}(d,c)_{n}}} & {\phantom{c}}
}}
\qquad
\vcenter{\xymatrix @C=3.5em{
d \ar[r]^-{\smash[t]{\beta_{0}\circ\alpha_{0}}} & c\,.
}}
$$
Note that $\beta_{j}\circ\alpha_{j}=\beta_{0}\circ\alpha_{0}$ for $0\leq j\leq n$, so that $j=0$
plays no prominent rôle.
\smallbreak
\item [(iii)] Define the more $\natural$ diagram $\EEDCpr{D}{C}=\EEDCpr{D}{C}(\qquest,?)$ without the op's, that is,
$$
\itemspace \EEDCpr{D}{C}:=B(\qquest\comma\cat{D}\commaC{C}?)\colon\;
\cat{D}\op\times\cat{C}\longrightarrow\sSets,\qquad(d,c)\longmapsto B(d\comma\cat{D}\commaC{C}c)\,,
$$
and define the morphism $\vthetaDCpr{D}{C}\colon\EEDCpr{D}{C}\longrightarrow
\FFDC{D}{C}$ similarly.
\end{itemize}
We write without decoration $\natural$ the diagram involving the op's since such diagrams are more commonly used, as we
shall see in Section~\ref{s-hocolim-colim-Q} for instance.
\end{Not}

\begin{Rem}
When $\cat{D}=\cat{C}$ in Notation~\ref{commaC-not}, for $c',c\in\cat{C}$, the category $c'\comma\cat{C}
\commaC{C}c$ is the usual double-comma category $c'\comma\cat{C}\comma c$. In this case, we write
$$
\EEC{C}:=\EEDC{C}{C}=B(\qquest\comma\cat{C}\comma{}?)\op
\qquad\mbox{and}\qquad
\vthetaC{C}:=\vthetaDC{C}{C}\,.
$$
We define $\EECpr{C}$ and $\vthetaCpr{C}$ similarly.
\end{Rem}

\begin{Rem}
\label{EC-interpret-rem}
For a subcategory $\cat{D}$ of a small category $\cat{C}$, the diagram $\EEDCpr{D}{C}
(\qquest,?)$ identifies canonically with the one given at level $n\in\bbN$ by
$$
\big\{\qquest\stackrel{\alpha_{0}}{\longrightarrow}d_{0}\stackrel{\gamma_{1}}{\longrightarrow}d_{1}
\stackrel{\gamma_{2}}{\longrightarrow}\ldots\stackrel{\gamma_{n}}{\longrightarrow}d_{n}
\stackrel{\beta_{n}}{\longrightarrow}?\big\}\,;
$$
we mean that all such morphisms $\alpha_{0}$ and $\gamma_{j}$ in $\cat{D}$ and $\beta_{n}$ in $\cat{C}$
are considered (with all possible objects $d_{j}$ in $\cat{D}$), and the face maps are given by composing two
successive maps and the degeneracies by inserting identities, but only among the $\gamma_{j}$'s.
\end{Rem}

\begin{Rem}
\label{cov-class-rem}
When $\cat{C}$ has an initial object $c_{0}$, the functor $\EECpr{C}(c_{0},?)$  is canonically isomorphic
to $B(\cat{C}\comma{}?)$ and its realization is sometimes called the ``covariant classifying $\cat{C}$-space''.
If $\cat{C}$ has a terminal object $\cterminal$, then $\EECpr{C}(\qquest,\cterminal)\cong B(\qquest\comma\cat{C})$ and
its realization is sometimes called the ``contravariant classifying $\cat{C}$-space''.
\end{Rem}

\begin{Thm}
\label{EC-cofibrant-thm}
Let $\cat{D}$ be a subcategory of a small category $\cat{C}$. Then, the morphism
$$
\vthetaDC{D}{C}\colon\EEDC{D}{C}\longrightarrow
\FFDC{D}{C}
$$
is a cofibrant approximation in the model category
$\UUU{\sSets}{\cat{D}\op\times\cat{C}}{\cat{D}\op\times\cat{D}}$. One has slightly more\,: $\EEDC{D}{C}$ is cofibrant
in $\UUU{\sSets}{\cat{D}\op\times\cat{C}}{\cat{D}\op\times\cat{D}}$ and $\vartheta$ is a weak equivalence $\cat{D}\op
\times\cat{C}$-objectwise. The same
holds for the morphism $\vthetaDCpr{D}{C}\colon\EEDCpr{D}{C}
\longrightarrow\FFDC{D}{C}$.
\end{Thm}

The proof of this theorem will follow the principle of \cite[\S\,9]{dugg} and will be presented
once a few technical lemmas are established. Combined with Theorem~\ref{cof-approx-USCD-thm},
Theorem~\ref{EC-cofibrant-thm} provides our main explicit cofibrant approximation
in $\USCD{S}{C}{D}$. We give it a name.

\begin{CorDef}
\label{main-cor-def}
Keeping notations as in Theorems \ref{cof-approx-USCD-thm} an~\ref{EC-cofibrant-thm}, for every diagram
$X\in\cat{S}^{\cat{C}}$, we denote by
$$
\xiCDXbar{C}{D}{X}\colon\QDCbar{D}{C}X:=\EEDC{D}{C}\udCDdot{C}{D}{\otimes}qX\longrightarrow X
$$
the corresponding functorial cofibrant approximation of $X$ in $\USCD{S}{C}{D}$. We call $\QDCbar{D}{C}X$
the \emph{bar cofibrant approximation of $X$} in $\USCD{S}{C}{D}$. Similarly, we write
$$
\xiCDXbarpr{C}{D}{X}\colon\QDCbarpr{D}{C}X:=\EEDCpr{D}{C}\udCDdot{C}{D}{\otimes}
qX\longrightarrow X
$$
and call it the \emph{opposite bar cofibrant approximation of $X$} in $\USCD{S}{C}{D}$. When $\cat{D}=
\cat{C}$, we also write $(\QCbar{C},\xiCbar{C})$ and $(\QCbarpr{C},\xiCbarpr{C})$. They indeed are
cofibrant approximations.
\qed
\end{CorDef}
\begin{Rem}
\label{unravel-rem}
Unravelling the definition of $\EEDC{D}{C}$ and that of $qX$, we see that
$$
\QDCbar{D}{C}X(?)=B(\qquest\comma\cat{D}\commaC{C}?)\op\udCDdot{C}{\smallqquest\in D}{\otimes}Q_{\cat{S}}
\big(X(\qquest)\big)\,,
$$
and similarly for $\QDCbarpr{D}{C}X(?)$, but without the op's.
\end{Rem}

\begin{Rem}
Since the realization of the nerve of a small category and of its opposite are homeomorphic,
for the model category of compactly generated Hausdorff spaces or that of spectra of such spaces
(with the strict or the stable model structure), the cofibrant approximations $\QDCbar{D}{C}$ and
$\QDCbarpr{D}{C}$ provide isomorphic diagrams.
\end{Rem}

\begin{Cor}
\label{quote-cor}
Keep notations as in Theorems \ref{cof-approx-USCD-thm} and~\ref{EC-cofibrant-thm}. Assume that $X\in\cat{S}^{\cat{C}}$
is \emph{$\cat{C}$-objectwise cofibrant}. Define the diagram $\nQX{X}\in\cat{S}^{\cat{C}}$ by $\nQX{X}:=\EEDC{D}{C}
\udCDdot{C}{D}{\otimes}\res_{\cat{D}}^{\cat{C}}X$, \ie
$$
\nQX{X}(?)=B(\qquest\comma\cat{D}\commaC{C}?)\op\udCDdot{C}{\smallqquest\in D}{\otimes} X(\qquest)
$$
and the morphism $\nzetaX{X}\colon\nQX{X}\too X$ by the composition of Corollary~\ref{X-objcof-cor}\,:
$$
\xymatrix @C=3.2em @R=1.6em{
{\EEDC{D}{C}}\,\udCDdot{C}{D}{\otimes}\res_{\cat{D}}^{\cat{C}}X
\ar[r]^-{\vartheta\smalludCDdot{C}{D}{\otimes}\res\nnspace X}
\ar@/_2em/[rrr]|{\ \nzetaX{X}\ }
& \FFDC{D}{C}\ \udCDdot{C}{D}{\otimes}\,
\res_{\cat{D}}^{\cat{C}}X
\ar[r]^-{\lambda_{\res\nnspace X}^{\cat{D},\cat{C}}}_-{\cong}
&\ind_{\cat{D}}^{\cat{C}}\res_{\cat{D}}^{\cat{C}}X
\ar[r]^-{\epsilon_{X}}
& X\,.
}
$$
Then, $\nzetaX{X}\colon\nQX{X}\too X$ is a cofibrant approximation of $X$ in $\USCD{S}{C}{D}$.
\end{Cor}

\proof
This follows from the above Corollary~\ref{X-objcof-cor} and from Theorem~\ref{EC-cofibrant-thm}.
\qed

\goodbreak \bigbreak
\centerline{*\ *\ *}
\medbreak

The following two lemmas constitute preparatory material for the proof of Theorem~\ref{EC-cofibrant-thm}.
For the first one, recall Notation~\ref{other-induced-not}.

\begin{Lem}
\label{useful-lem}
Let $\cat{D}$ be a full subcategory of a small category $\cat{C}$. We have\,:
\begin{itemize}
\item [(i)] Let $\cat{S}$ be a cofibrantly generated simplicial model category.
If $X$ is a cofibrant object in the model category $\USCD{S}{C}{D}$, then the morphism
$$
\incl\uCdot{\pair{\;\,}{C}}{\odot}X
\colon\quad\partial\Delta^{n}\uCdot{\pair{\;\,}{C}}{\odot}X
\longrightarrow\Delta^{n}\uCdot{\pair{\;\,}{C}}{\odot} X
$$
is a cofibration in $\USCD{S}{C}{D}$ for each $n\in\bbN$.
\smallbreak
\item [(ii)]  Let $\{d_{i}\}_{i\in I}$ a set
of objects of $\cat{D}$. Then the diagram
$$
\itemspace\coprod_{i\in I}\ {\FFC{C}(\qquest,d_{i})}\in\ \sSets^{\cat{C}\op}
$$
is a cofibrant object in the model category $\USCD{{}\sSets}{C\op}{D\op}$.
\smallbreak
\item [(iii)] In a model category, an object $X$ which is a sequential colimit
$$
\itemspace X=\colim_{n\in\bbN}\,\big(X_{0}\longrightarrow\ldots\longrightarrow X_{n}\longrightarrow X_{n+1}
\longrightarrow\ldots\big)
$$
of cofibrations with $X_{0}$ cofibrant is itself cofibrant.
\end{itemize}
\end{Lem}

\Prf
(i) As explained in Theorem~\ref{USCD-simpl-thm}, the model category $\USCD{S}{C}{D}$ inherits
a canonical structure of simplicial model category with ``action'' given by $\uCdot{\pair{\;\,}
{C}}{\odot}$. So, the result now follows from Remark~\ref{flat-rem} or~\cite[Prop.\ 9.3.9\,(1)\,(a)]{hirsch} and
the fact that the inclusion map $\incl\colon\partial\Delta^{n}\longrightarrow\Delta^{n}$ is a cofibration
of simplicial sets.

(ii) follows from Lemma \ref{mor-objwise-cof-lem} and the general fact that, in a model category,
the class of cofibrations is determined by the left lifting property with respect to some fixed class of morphisms
(trivial fibrations).
This forces the coproduct, the push-out and the sequential colimit of cofibrations to be again a cofibration. This
gives us (iii) as well (adding the cofibration $\varnothing\to X_{0}$ at the beginning, if one prefers).
\qed

\begin{Lem}
\label{skeleton-lem}
Let $K$ be a simplicial set. For $n\in\bbN$, denote by $\nd_{n}(K)$ the set of non-degenerate $n$-simplices
of $K$. Then, there is a canonical isomorphism
$$
K\cong\colim_{n\in\bbN}\,\sk_{n}(K)
$$
of simplicial sets, where $\sk_{0}(K):=\Delta^{0}\times\nd_{0}(K)=\underline{\nd_{0}(K)}=\underline{K_{0}}$
and for $n\geq 1$, $\sk_{n}(K)$ is determined by a push-out
$$
\xymatrix{
\partial\Delta^{n}\times\nd_{n}(K) \ar[r] \ar[d] & \sk_{n-1}(K) \ar@{-->}[d] \\
\Delta^{n}\times\nd_{n}(K) \ar@{-->}[r] & \sk_{n}(K)
}
$$
and the ``structure map'' $\sk_{n-1}(K)\longrightarrow\sk_{n}(K)$ for the colimit is given by the right
vertical map in this push-out.
\end{Lem}

Note that $|\sk_{n}(K)|$ is the $n$-skeleton of $|K|$ (with its canonical CW-structure), hence
the notation\,: see Goerss-Jardine~\cite[p.\ 8]{goja}, where a proof of~\ref{skeleton-lem} can be found.

\smallbreak

\PrfOf{Theorem \ref{EC-cofibrant-thm}}
Let us first show that $\vthetaDCpr{D}{C}$ is a weak equivalence. Fix two objects
$d$ in $\cat{D}$ and $c$ in $\cat{C}$. We have to show that the map of simplicial sets
$$
\vthetaDCpr{D}{C}(d,c)\colon B(d\comma\cat{D}\commaC{C}c)\longrightarrow
\underline{\FDC{D}{C}(d,c)}
$$
is a weak equivalence. Consider $\FDC{D}{C}(d,c)$ as a discrete category.
Plainly, we have a canonical isomorphism of simplicial sets
$$
B\FDC{D}{C}(d,c)\cong\underline{\FDC{D}{C}(d,c)}\,,
$$
considered below as an identity and written ``$=$''. Let
$$
\pi\colon d\comma\cat{D}\commaC{C}c\longrightarrow\FDC{D}{C}(d,c),\qquad
\big(d\stackrel{\alpha}{\rightarrow}d'\stackrel{\beta}{\rightarrow}c\big)
\longmapsto\big(d\stackrel{\beta\circ\alpha}{\rightarrow}c\big)
$$
be the ``composition functor'' and let the ``pre-insert identity'' functor be
$$
\iota\colon\FDC{D}{C}(d,c)\longrightarrow d\comma\cat{D}\commaC{C}c,\qquad
\big(d\stackrel{\gamma}{\rightarrow}c\big)\longmapsto\big(d\stackrel{\;\id_{d}}{\rightarrow}d
\stackrel{\gamma}{\rightarrow}c\big)\,.
$$
Clearly, the equalities $\pi\circ\iota=\id_{\mor_{\cat{C}}(d,c)}$ and $B\pi=\vthetaDCpr{D}{C}(d,c)$ hold.
Now, we define a natural transformation $\nu\colon\iota\circ\pi\longrightarrow\id_{(d\smallcomma\cat{D}\smallcommaC{C}c)}$
by the assignment
$$
\vcenter{\xymatrix{
d \ar[r]^-{\alpha} & d' \ar[r]^-{\beta} & c
}}
\qquad
\vcenter{\xymatrix{
{} \ar@{|->}[r]^{\nu} & {}
}}
\qquad
\vcenter{\xymatrix @R=1.5em{
d \ar[r]^-{\id_{d}} \ar[d]_{\id_{d}\!} & d \ar[r]^-{\beta\circ\alpha} \ar[d]_{\alpha\!} & c
\ar[d]^{\!\id_{c}} \\
d \ar[r]^-{\alpha} & d' \ar[r]^-{\beta} & c\,.
}}
$$
We deduce that $B\pi\circ B\iota=B(\pi\circ\iota)=B\id_{\mor_{\cat{C}}(d,c)}=\id_{\underline
{\mor_{\cat{C}}(d,c)}}$ and, by Segal's argument \cite[Prop.\ 2.1]{segal}, $\nu$ induces a homotopy equivalence
$$
B\iota\circ B\pi=B(\iota\circ\pi)\sim B\id_{(d\smallcomma\cat{D}\smallcommaC{C}c)}=\id_{B(d\smallcomma\cat{D}\smallcommaC{C}c)}\,,
$$
showing that $B\iota$ and $B\pi=\vthetaDCpr{D}{C}(d,c)$
are mutual inverse homotopy equivalences,
hence weak equivalences.

\medbreak

Now, we prove that $\EEDCpr{D}{C}$ is a cofibrant object in $\USCD{{}\sSets}{D\op\times\cat{C}}{D\op\times\cat{D}}$. We
apply Lemma \ref{skeleton-lem} $\cat{D}\op\!\times\!\cat{C}$-objectwise,
namely, for each pair $(d,c)\in\cat{D}\op\times\cat{C}$, we deduce from this Lemma that we have an isomorphism
$$
\EEDCpr{D}{C}(d,c)=B(d\comma\cat{D}\commaC{C}c)\cong\colim_{n\in\bbN}\,\sk_{n}\big(B(d\comma\cat{D}\commaC{C}c)\big)\,,
$$
of simplicial sets, where
$$
\sk_{0}\big(B(d\comma\cat{D}\commaC{C}c)\big)=\Delta^{0}\times\nd_{0}\big(B(d\comma\cat{D}\commaC{C}c)\big)
=\underline{\obj(d\comma\cat{D}\commaC{C}c)}
$$
and where, for each $n\geq 1$, $\sk_{n}$ is obtained from $\sk_{n-1}$ by a push-out
$$
\xymatrix{
\partial\Delta^{n}\times\nd_{n}\big(B(d\comma\cat{D}\commaC{C}c)\big) \ar[r] \ar[d]
& \sk_{n-1}\big(B(d\comma\cat{D}\commaC{C}c)\big) \ar[d]
\\
\Delta^{n}\times\nd_{n}\big(B(d\comma\cat{D}\commaC{C}c)\big) \ar[r] & \sk_{n}\big(B(d\comma\cat
{D}\commaC{C}c)\big)
}
$$
Observe from Remark \ref{EC-interpret-rem} that for each $n\in\bbN$, one has a canonical bijection
$$
\nd_{n}\big(B(d\comma\cat{D}\commaC{C}c)\big)
\cong
\coprod_{(d_{0}\rightarrow\ldots\rightarrow d_{n})\in S_{n}}
\mor_{\cat{D}\times\cat{C}\op}\!\big((d,c),(d_{0},d_{n})\big)
$$
where $S_{n}$ designates the set
$$
\big\{d_{0}\stackrel{\gamma_{1}}{\longrightarrow}d_{1}\stackrel{\gamma_{2}}{\longrightarrow}\ldots
\stackrel{\gamma_{n}}{\longrightarrow}d_{n}\,\big|\,d_{j}\in\obj(\cat{D})\;\,\mbox{and}\;\,\gamma_{i}\in
\arr(\cat{D})\;\,\mbox{with}\;\,\gamma_{i}\neq\id_{d_{i}}\big\}\,.
$$
Observe moreover that for a simplicial set $K$ and a set $S$, one has
$$
K\times S=K\times\underline{S}=K\odot\underline{S}\,,
$$
where $\odot$ is as in Example \ref{coupling-ex} (with $\cat{S}$ standing for $\sSets$) and $\underline{S}$
is ``simplicially constant''. So our skeleton decomposition of $\EEDCpr{D}{C}(d,c)=B(d\comma\cat{D}\commaC{C}c)$
becomes
$$
\sk_{0}\big(\EEDCpr{D}{C}(d,c)\big):=
\coprod_{d_{0}\in\cat{D}}\underline{\mor_{\cat{D}\times\cat{C}\op}\!\big((d,c),(d_{0},d_{0})\big)}
$$
and for each $n\geq 1$, we have a push-out square
$$
\xymatrix{
\partial\Delta^{n}\ \odot
{\displaystyle\coprod_{\smash{(d_{0}\rightarrow\ldots\rightarrow d_{n})\in S_{n}}}^{\vphantom{s}}}
\underline{\mor_{\cat{D}\times\cat{C}\op}\!\big((d,c),(d_{0},d_{n})\big)}
\ar[r] \ar[d]
&   \sk_{n-1}\big(\EEDCpr{D}{C}(d,c)\big) \ar[d]
\\
\Delta^{n}\ \odot
{\displaystyle\coprod_{\smash{(d_{0}\rightarrow\ldots\rightarrow d_{n})\in S_{n}}}^{\vphantom{s}}}
\underline{\mor_{\cat{D}\times\cat{C}\op}\!\big((d,c),(d_{0},d_{n})\big)}
\ar[r]
& \sk_{n}\big(\EEDCpr{D}{C}(d,c)\big)\,.
}
$$

Let us rewrite the above in a more diagram-category language, using Notation~\ref{mor-eta-not}. We have proven that our
diagram $\EEDC{D}{C}(\qquest,?)\in\sSets^{\cat{D}\op\times\cat{C}}$ is isomorphic to a sequential colimit
$$
\EEDC{D}{C}\cong\colim_{n\in\bbN}\,\mathcal{X}_{n}
$$
of diagrams $\mathcal{X}_{n}\in\sSets^{\cat{D}\op\times\cat{C}}$, where $\mathcal{X}_{n}(\qquest,?):=\sk_{n}\big(B
(\qquest\comma\cat{D}\commaC{C}?)\big)$ and we have
$$
\mathcal{X}_{0}(\qquest,?)=
\coprod_{d_{0}\in\cat{D}}\FFC{{}(\cat{D}\times\cat{C}\op)}\!\big((\qquest,?),(d_{0},d_{0})\big)
$$
and, for each $n\geq 1$, we have a push-out square (we call the left-hand map $j_{n}$)
$$
\xymatrix@R=2em{
\partial\Delta^{n}\hspace*{-1.6em}
\uCdot{\pair{\quad\quad}{\cat{D}\op\times\cat{C}}}{\odot}
{\displaystyle\coprod_{\smash{(d_{0}\rightarrow\ldots\rightarrow d_{n})\in S_{n}}}^{\vphantom{s}}}
\FFC{{}(\cat{D}\times\cat{C}\op)}\!\big((\qquest,?),(d_{0},d_{n})\big)
\ar[r] \ar[d]^{=:\displaystyle j_{n}}
& {\mathcal{X}}_{n-1}(\qquest,?)\ar[d]
\\
\Delta^{n}\hspace*{-1.6em}\uCdot{\pair{\quad\quad}{\cat{D}\op\times\cat{C}}}{\odot}
{\displaystyle\coprod_{\smash{(d_{0}\rightarrow\ldots\rightarrow d_{n})\in S_{n}}}^{\vphantom{s}}}
\FFC{{}(\cat{D}\times\cat{C}\op)}\!\big((\qquest,?),(d_{0},d_{n})\big)
\ar[r]
& {\mathcal{X}}_{n}(\qquest,?)\,.
}
$$
We want to conclude by means of Lemma~\ref{useful-lem}\,(iii) and we have to check two things\,:
\begin{itemize}
\item[(a)] the first object $\mathcal{X}_{0}\in\sSets^{D\op\times\cat{C}}$ is cofibrant in $\USCD{{}\sSets}{D\op\times
\cat{C}}{D\op\times\cat{D}}$;
\item[(b)] $\mathcal{X}_{n-1}\too\mathcal{X}_{n}$ is a cofibration in $\USCD{{}\sSets}{D\op\times\cat{C}}{D\op\times
\cat{D}}$ for all $n\geq1$.
\end{itemize}
Both are taken care of by Lemma~\ref{useful-lem}.
Property (a) follows from its part (ii) applied to the pair $\cat{D}\times\cat{D}\op\ \subset\ \cat{D}\times\cat{C}\op$
instead of $\cat{D}\subset\cat{C}$. The same argument guarantees that the coproduct in the above push-out square is a
cofibrant object in  $\USCD{{}\sSets}{D\op\times\cat{C}}{D\op\times\cat{D}}$, which is a simplicial model category by
Theorem \ref{USCD-simpl-thm}. So, the morphism $j_{n}$ is a cofibration by part~(i) of Lemma~\ref{useful-lem} and so is
$\mathcal{X}_{n-1}\too\mathcal{X}_{n}$ by push-out, hence (b) above. This gives the result for $\vthetaDCpr{D}{C}\colon
\EEDCpr{D}{C}\too\FFDC{D}{C}$.

\smallbreak

The proof for $\EEDC{D}{C}$ and $\vthetaDC{D}{C}$ is similar, {\sl mutatis} op's {\sl mutandis}.
\qed

\goodbreak \bigbreak
\centerline{*\ *\ *}
\smallbreak

\begin{Rem}
\label{various-rem}
Let $\cat{D}$ be a full subcategory of $\cat{C}$.
With Remark~\ref{cof-approx-USCD-rem} in mind, the reader might ask whether one obtains a new cofibrant approximation
by the method exposed there of inducing up an approximation $\bbE_{0}\to\FFC{D}$ from the absolute situation
$\USC{{}\sSets}{D\op\times\cat{D}}$ to an approximation in the relative one $\USCD{{}\sSets}{D\op\times\cat{C}}
{D\op\times\cat{D}}$, when applied to the above $\EEDC{D}{C}$. We leave it to the reader to see that nothing new
happens in this way, that is, the induction of $\EEDC{D}{D}$ is isomorphic to $\EEDC{D}{C}$ in a compatible
way with the $\vthetaC{{}}$'s, and similarly for $\EECpr{{}}$ and $\vthetaCpr{{}}$.

One can also wonder if other (absolute) approximations of $\FFC{D}$, for instance the two ones of Dugger~\cite{dugg},
will produce essentially different approximations after induction to the relative model category. Strictly speaking
the answer is yes, although inducing up the one of \cite[Lem.\,2.7]{dugg}, for instance, only differs from the above
$\EEDC{D}{C}$ by an ``edgewise subdivision''. The details are again left to the interested reader.
\end{Rem}


\section{Homotopy colimits versus colimits of cofibrant approximations}

\label{s-hocolim-colim-Q}


%
Fix a cofibrantly generated simplicial model category $\cat{S}$ and a small category~$\cat{C}$. We compare the two
possible approaches to homotopy colimits. The reader should have in mind the identification $*\,\otimes_{\cat{C}}-
\cong\colim_{\cat{C}}(-)$ of Example~\ref{colim=*tensor(-)-ex}, where $*\in\sSets^{\cat{C}\op}$ is the constant
diagram taking the point $\Delta^{0}\in\sSets$ as value. Comma categories are defined in~\ref{commaC-not}; for a
diagram $X\in\cat{S}^{\cat{C}}$, the objectwise cofibrant approximation $qX$ and the cofibrant approximation
$\QCbar{C}X$ are introduced in \ref{mor-eta-not} and \ref{main-cor-def} respectively.

\begin{Def}
\label{hocolim-Lcolim-def}
Consider a diagram $X\in\cat{S}^{\cat{C}}$. Let us fix the terminology.
\begin{itemize}
\item [(i)] The \emph{homotopy colimit} of $X$ is the object of $\cat{S}$ given by
$$
\itemspace\hocolim_{\cat{C}}X:=B(\qquest\comma\cat{C})\op\,\dCdot{\qquest\in C}{\otimes}X(\qquest)\,.
$$
\item [(ii)] The \emph{(bar) $L$-colimit} of $X$ is the object of $\cat{S}$ given by
$$
\itemspace\Lcolim_{\cat{C}}X:=\colim_{\cat{C}}\QCbar{C}X
\,.
$$
\end{itemize}
Both constructions are clearly functorial. There are natural
morphisms (in $\cat{S}$)
$$
\xymatrix @C=3.5em{
\hocolim_{\cat{C}}X=B(\qquest\comma\cat{C})\op{\smash[b]{\dCdot{\qquest\in C}{\otimes}X(\qquest)}} \ar[r]^-{\Projec
{\otimes}X} & {*}\,{\smash[b]{\dCdot{C}{\otimes}}}X=\colim_{\cat{C}}X\,,
}
$$
where $\Projec\colon B(\qquest\comma\cat{C})\op\too *$ is the obvious $\cat{C}\op$-objectwise constant map, and
$$
\xymatrix@C=4.5em{
\Lcolim_{\cat{C}}X
\ar[r]^-{\colim_{\cat{C}}\xiCXbar{C}{X}} &
\colim_{\cat{C}}X\,,
}
$$
where $\xiCXbar{C}{X}\colon\QCbar{C}X\too X$ is the morphism of~\ref{main-cor-def}, see also Corollary~\ref{cof-approx-USC-cor}.
\end{Def}

\begin{Thm}
\label{Lcolim=hocolim-thm}
For every diagram $X\in\cat{S}^{\cat{C}}$, there is a commutative diagram
$$
\xymatrix@C=4.2em{
\Lcolim_{\cat{C}}X \ar@{=}[r]^-{\sim} \ar[rd]_{\colim_{\cat{C}}\xiCXbar{C}{X}\,} & \hocolim_{\cat{C}}qX
\ar[r]^-{\hocolim\eta_{X}}
& \hocolim_{\cat{C}}X
\ar[ld]^{\Projec{\otimes}X} \\
& \colim_{\cat{C}}X & \\
}
$$
Moreover, the following properties are equivalent\,:
\begin{itemize}
\item[(a)]
The horizontal arrow $\Lcolim_{\cat{C}}X\too \hocolim_{\cat{C}}X$ is a natural weak equivalence in $\cat{S}$ for
all $X\in\cat{S}^{\cat{C}}$.
\item[(b)] The functor $\hocolim_{\cat{C}}\colon\cat{S}^{\cat{C}}\too\cat{S}$ is weakly homotopy invariant, that is,
takes $\cat{C}$-\weqs\ to \weqs.
\end{itemize}
\end{Thm}

\Prf
Let us start with the following two observations\,:
\begin{itemize}
    \item [(1)] The tensor product commutes with colimits, as can be checked directly from
    the definition or by adjunction; in particular, for every $\mathcal{X}\in\cat{S}^{\cat{C}\op\times\cat{C}}$ and
    $Y\in\cat{S}^{\cat{C}}$,
    we have a natural isomorphism $(\colim_{\cat{C}}\mathcal{X})\dCdot{C}{\otimes}Y\cong
    \colim_{\cat{C}}(\mathcal{X}\udCDdot{C}{C}{\otimes}Y)$.
    \item [(2)] We have an isomorphism $\colim_{?\in\cat{C}}B(\qquest\comma\cat{C}\comma{}?)\op\cong B(\qquest\comma\cat{C})\op$.
    This can be checked directly from the universal property of the colimit.
\end{itemize}
We compute
$$
\arraycolsep1pt
\renewcommand{\arraystretch}{1.5}
\begin{array}{rclll}
\Lcolim_{\cat{C}}X & =
& \colim_{?\in\cat{C}}\big(\QCbar{C}X(?)\big)& \qquad & \mbox{(by definition)}
\\
& = & \colim_{?\in\cat{C}}\Big(B(\qquest\comma\cat{C}\comma{}?)\op\,\udCDdot{C}{{}\qquest\in\cat{C}}{\otimes}qX(\qquest)\Big)
& \qquad & \mbox{(by definition of } \QCbar{C}X\,)
\\
& \cong & \big(\colim_{?\in\cat{C}}B(\qquest\comma\cat{C}\comma{}?)\op\big)\,\dCdot{{}\qquest\in\cat{C}}{\otimes}qX(\qquest) & &
\mbox{(by~(1) above)}
\\
& \cong & B(\qquest\comma\cat{C})\op\,\dCdot{{}\qquest\in\cat{C}}{\otimes}qX(\qquest) & & \mbox{(by~(2) above)}
\\
& = & \hocolim_{\cat{C}}qX & & \mbox{(by definition).} \\
\end{array}
$$
Hence the isomorphism $\Lcolim_{\cat{C}}X\cong \hocolim_{\cat{C}}qX$ of the statement. The commutativity of the diagram
is left as an exercise.
Let us see that (a) and (b) are equivalent. This is now easy. Since $\Lcolim_{\cat{C}}$ is a left derived functor
it preserves weak-equivalences, hence (a)$\Longrightarrow$(b). Conversely, since $\eta_{X}$ is a weak equivalence
in $\cat{S}^{\cat{C}}$, one clearly deduce (b)$\Longrightarrow$(a) from the above.
\qed

\medskip

Let us mention some important cases where the theorem applies.

\begin{Prop}
For the following (cofibrantly generated) simplicial model categories $\cat{S}$, the functor $\hocolim_
{\cat{C}}\colon\USC{S}{C}\longrightarrow\cat{S}$ is weakly homotopy invariant, for every small category
$\cat{C}$\,:
\begin{itemize}
\item [(i)] the model category of simplicial sets;
\item [(ii)] the model category of pointed simplicial sets;
\item [(iii)] the strict model category of spectra of simplicial sets;
\item [(iv)] the stable model category of spectra of simplicial sets.
\end{itemize}
\end{Prop}

\Prf
This is well-known to the experts. Noting that every simplicial set and every pointed simplicial
set is cofibrant (cf.\ \cite[Prop.\ 3.2.2 and Cor.\ 3.6.6]{hov}), the result for $\sSets$ and for
pointed simplicial sets follows from \cite[Thm.\ 19.4.2 (1)]{hirsch} (a proof for $\sSets$ is already
contained in \cite[Lem., p.\ 329]{boukan}). For the strict model category of spectra, this follows from
the case of pointed simplicial sets, since strict weak equivalences are defined levelwise and since homotopy
colimits can be taken levelwise, as easily checked. For the stable model category of spectra, this
is the content of \cite[Lem.\ 5.18]{thom}.
\qed

\begin{Rem}
\label{zig-zag-rem}
If we had defined $\Lcolim$ using the diagram $B(\qquest\comma\cat{C}\comma{}?)$ in place of
$B(\qquest\comma\cat{C}\comma{}?)\op$, we would get a non-canonical
zig-zag of weak equivalences instead of the canonical isomorphism in Theorem~\ref{Lcolim=hocolim-thm}, still with a
commutative diagram as in
the statement. More generally, since two left derived functors of the same functor are weakly equivalent, one would
have such a zig-zag of weak-equivalences even if we replace $\QCbar{C}X$ by an arbitrary functorial cofibrant approximation of $X$.
\end{Rem}


\appendix


\section{Simplicial model categories}

\label{AppIIA}


%
We define simplicial model categories; examples of such categories of ``spaces'' and
of spectra will be presented in Appendix A of \cite{bamaIII} with some details.

\medskip

Let $\sSets$ denote the category of simplicial sets. Given two simplicial sets $K$ and $K'$, we let
$\Map_{\sSets}(K,K')=\Map_{\sSets} (K,K')_{\bullet}$ be the simplicial set defined by
$$
\Map_{\sSets}(K,K')_{q}:=\mor_{\sSets}(K\times\Delta^{\!q},K')\qquad(q\in\bbN)
$$
and with the obvious face and degeneracy maps coming from the simplicial structure of
$\Delta^{\!q}=\mor_{\DDelta\op}(-,[q])$.

\begin{Def}
\label{simpl-mod-cat-def}
A \emph{simplicial model category} is a model category $\cat{M}$ equipped with
\begin{itemize}
\item [(i)] a simplicial set $\Map(m_{1},m_{2})$, for all $m_{1},m_{2}\in\cat{M}$; \item [(ii)] a
``composition'' map of simplicial sets
$$
\itemspace\Map(m_{2},m_{3})\times\Map(m_{1},m_{2})\longrightarrow\Map(m_{1},m_{3})\,,
$$
for all $m_{1},m_{2},m_{3}\in\cat{M}$; \item [(iii)] a bijection
$\mor_{\cat{M}}(m_{1},m_{2})\cong\Map(m_{1},m_{2})_{0}$ (degree-zero part), that is compatible with
the compositions, for all $m_{1},m_{2}\in\cat{M}$;
\end{itemize}

moreover, the following axioms are required\,:
\begin{itemize}
\item [(1)] the ``composition'' map is associative with two sided identity given, for each $m\in\cat{M}$, by the
image of $\id_{m}\in\mor_{\cat{M}}(m,m)$ in $\Map(m,m)$;
\item [(2)] for $K\in\sSets$ and $m\in\cat{M}$, there
exists an object $K\odot m\in\cat{M}$ such that, for all $\ell\in\cat{M}$, there is an isomorphism
$$
\itemspace\Map(K\odot m, \ell)\cong\Map_{\sSets}(K,\Map(m,\ell))
$$
of simplicial sets, natural in $\ell$;
\item [(3)] for $K\in\sSets$ and $m\in\cat{M}$, there exists an object
$m^{K}\in\cat{M}$ such that, for all $\ell\in\cat{M}$, there is an isomorphism
$$
\Map(K\odot \ell, m)\cong\Map(\ell,m^{K})
$$
of simplicial sets, natural in $\ell$;
\item [(4)] for every cofibration $p\colon\ell\longrightarrow\ell'$ in
$\cat{M}$ and for every fibration $i\colon m\longrightarrow m'$ in $\cat{M}$, the induced map of simplicial sets
$$
\itemspace
\Map(\ell',m)\longrightarrow\Map(\ell,m)\times_{\Map(\ell,m')}\Map(\ell',m')
$$
is a fibration; it is a trivial fibration if either $i$ or $p$ is moreover trivial, \ie also a weak equivalence.
\end{itemize}
\end{Def}

In this definition, we followed Hirschhorn \cite[\S\,9.1]{hirsch}, but \emph{not} Quillen \cite{quil},
where he only requires axioms (2) and (3) for $K$ a \emph{finite} simplicial set; as a consequence, the category
$\Top$ of all topological spaces is, for us, no simplicial model category.

\begin{Rem}
\label{odot-nat}
For a simplicial model category $\cat{M}$, one shows that the assignment
$$
\sSets\times\cat{M}\stackrel{\odot}{\longrightarrow}\cat{M},\qquad
(K,m)\longmapsto K\odot m\,,
$$
given by axiom (2) in the latter definition, is a (bi-)functor; we will refer to it as the ``action'' of $\sSets$
on $\cat{M}$. Compare Example~\ref{coupling-ex} above.
\end{Rem}

\begin{Ex}
The model category of simplicial sets is a simplicial model category; the corresponding
``action'' reads $K\odot Y=K\times Y$, for $K,Y\in\sSets$. Similarly, for the category of pointed
simplicial sets, we have $K\odot Y=K_{+}\wedge Y$. For the category of unpointed (resp.\ pointed)
compactly generated Hausdorff topological spaces, the ``action'' is given by $K\odot Y=|K|\times Y$
(resp.\ $K\odot Y=|K|_{+}\wedge Y$).
\end{Ex}

There is an obvious asymmetry in the above Definition~\ref{simpl-mod-cat-def}\, since we require the functor
$(\ell,m)\mapsto\Map(\ell,m)$ to have what we call in Definition~\ref{corner-map-def} the \emph{corner-map property}
but we do not require the other adjoint $(K,m)\mapsto m^{K}$ to have such a property. In fact both are equivalent
(see Hirschhorn~\cite[Prop.\,9.3.7]{hirsch})\,:

\begin{Prop}
\label{bal-corner-map-prop}
A model category satisfying (i)--(iii) and (1)--(3) of Definition~\ref{simpl-mod-cat-def} satisfies condition (4)
if and only if it satisfies the following condition\,:
\begin{itemize}
\item [($4^{\prime}$)] For every cofibration $i\colon K\to K'$ in $\sSets$ and for every fibration $p\colon m\to
\ell$ in $\cat{M}$ the induced morphism in $\cat{M}$
$$
\itemspace
m^{K'}\longrightarrow m^{K}\times_{\ell^{K}}\ell^{K'}
$$
is a fibration; it is a trivial fibration if either $i$ or $p$ is moreover trivial.
\end{itemize}
\end{Prop}

We use this flexibility to prove the following relative version of~\cite[Thm.\ 11.7.3]{hirsch}.

\begin{Thm}
\label{USCD-simpl-thm}
Let $\cat{M}$ be a cofibrantly generated simplicial model category and let
$
\odot\colon\sSets\times\cat{M}\longrightarrow\cat{M}
$
be the corresponding ``action'' as in~\ref{odot-nat}. Let $\cat{D}\subset\cat{C}$ be a pair of small categories.
Then the relative model structure $\USCD{M}{C}{D}$ on $\cat{M}^{\cat{C}}$ is a simplicial model category where the
``action'' is given by
$$
\uCdot{\pair{\;\,}{C}}{\odot}\colon\sSets\times\cat{M}^{\cat{C}}\longrightarrow\cat{M}^{\cat{C}},\qquad
(K,X)\longmapsto (K\uCdot{\pair{\;\,}{C}}{\odot}X)(?)=K\odot \big(X(?)\big)
$$
and where for every $K\in\sSets$ and every $X,Y\in\cat{M}^{\cat{C}}$, the object $X^{K}\in\cat{M}^{\cat{C}}$ and
the simplicial set $\Map(X,Y)\colon\DDelta\op\too\Sets$ are defined by
$$
X^{K}(?):=\big(\,X(?)\,\big)^{K}
\qquad
\mbox{and}
\qquad
\Map(X,Y)_{n}=\mor_{\cat{M}^{\cat{C}}}(\Delta^{n}\odot X,Y).
$$

\end{Thm}

\Prf
Hirschhorn establishes all the required adjunctions in~\cite[Thm.\ 11.7.3]{hirsch} and we are left to prove
condition~($4^{\prime}$) of Proposition~\ref{bal-corner-map-prop} for the relative structure, which goes as
in {\sl loc.\ cit}. Let $i\colon K\to K'$ be a cofibration in the model category $\sSets$ and let $p\colon X\to Y$ in
$\cat{M}^{\cat{C}}$
be a fibration $\cat{D}$-objectwise. We have to check that the induced map $\Phi\colon X^{K'}\longrightarrow X^{K}
\times_{Y^{K}}Y^{K'}$ is a fibration in $\USCD{M}{C}{D}$. For any $d\in\cat{D}$, $\Phi(d)=\varphi$ where $\varphi
\colon X^{K'}(d)=X(d)^{K'}\longrightarrow X(d)^{K}\times_{Y(d)^{K}}Y(d)^{K'}=(X^{K}\times_{Y^{K}}Y^{K'})(d)$ is
the corner map induced by $i\colon K\too K'$ and by $p(d)\colon X(d)\to Y(d)$, which is a fibration by assumption.
The result follows from the above condition~($4^{\prime}$) for $\cat{M}$. The same argument works if $i$ or $p$ is trivial.
\qed



\end{document}